\documentclass[12pt,amstex]{amsart}
       \usepackage{graphicx}
\usepackage{amssymb,amsthm,amsmath,epic,eepic,ulem, fullpage}
\usepackage{alltt}
\usepackage{setspace}

\newcommand{\tld}{\widetilde }
\theoremstyle{plain}
\newtheorem{thm}{Theorem}[section]

\newtheorem{cor}[thm]{Corollary}
\newtheorem{prop}[thm]{Proposition}

\newtheorem{lem}[thm]{Lemma}
\theoremstyle{definition}
\newtheorem{defn}[thm]{Definition}
\newtheorem{eg}[thm]{Example}
\theoremstyle{remark}
\newtheorem{rem}[thm]{Remark}

\newtheorem*{cl}{Claim}

\def\a{{\mathfrak{a}}} \def\b{{\mathfrak{b}}}

\def\F{{\mathbb{F}}}

\def\m{{\mathfrak{m}}}

\def\Z{{\mathbb{Z}}}\def\N{{\mathbb{N}}} \def\O{{\mathcal{O}}}
\def\Q{{\mathbb{Q}}} \def\R{{\mathbb{R}}} 
\newcommand{\DuBois}[1]{{\uuline \Omega {}^0_{#1}}}

 \def\Ker{{\mathrm{Ker}}} 
  \def\Ann{{\mathrm{Ann}}}
\def\Spec{{\mathrm{Spec\; }}} 
 \def\adj{{\mathrm{adj}}}
\DeclareFontFamily{OMS}{rsfs}{\skewchar\font'60}
\DeclareFontShape{OMS}{rsfs}{m}{n}{<-5>rsfs5 <5-7>rsfs7 <7->rsfs10 }{}
\DeclareSymbolFont{rsfs}{OMS}{rsfs}{m}{n}
\DeclareSymbolFontAlphabet{\scr}{rsfs}

\newcommand{\bQ}{\mathbb{Q}}
\newcommand{\bH}{\mathbb{H}}
\newcommand{\ba}{\mathfrak{a}}
\DeclareMathOperator{\coherent}{{coh}}
\DeclareMathOperator{\quasicoherent}{{qcoh}}
\newcommand{\mydot}{{{\,\begin{picture}(1,1)(-1,-2)\circle*{2}\end{picture}\ }}}
\newcommand{\qis}{\simeq_{\text{qis}}}
\newcommand{\sF}{\scr{F}}
\newcommand{\sG}{\scr{G}}
\newcommand{\sH}{\scr{H}}
\newcommand{\myR}{{\bf R}}
\DeclareMathOperator{\sHom}{{\sH}om}

\DeclareMathOperator{\exc}{{exc}}
\DeclareMathOperator{\Supp}{{Supp}}
\newcommand{\tensor}{\otimes}
\newcommand{\mJ}{\mathcal{J}}
\DeclareMathOperator{\red}{red}
\newcommand{\bC}{\mathbb{C}}
\newcommand{\bm}{\mathfrak{m}}
\newcommand{\bb}{\mathfrak{b}}

\newcommand{\Cech}{{$\check{\text{C}}$ech} }

\input xy
\xyoption{all}
\renewcommand{\DuBois}[1]{{\underline \Omega {}^0_{#1}}}

\begin{document}

\title{Rational singularities associated to pairs}
\author{Karl Schwede and Shunsuke Takagi}

\thanks{\\The first author was partially supported by RTG grant number 0502170 and
also as an National Science Foundation post-doc.\\
The second author was partially supported by Grant-in-Aid for Young
Scientists (B) 17740021 from JSPS and by Program for Improvement of Research
Environment for Young Researchers from SCF commissioned  by MEXT of Japan.}
\dedicatory{Dedicated to Professor Mel Hochster on the occasion of his sixty-fifth birthday.}
\maketitle

\section{Introduction and background}
\label{SectionIntroduction}

Rational singularities are a class of singularities which have been heavily studied since their introduction in the 1960s.  Roughly speaking, an algebraic variety has rational singularities if its structure sheaf has the same cohomology as the structure sheaf of a resolution of singularities.  Rational singularities enjoy many useful properties, in particular they are both normal and Cohen-Macaulay.  Furthermore, many common varieties have rational singularities, including toric varieties and quotient varieties.  Rational singularities are also known to be closely related to the singularities of the minimal model program.  In particular, it is known that log terminal singularities are rational and that Gorenstein rational singularities are canonical.

There is, however, an important distinction between rational singularities and singularities of the minimal model program.  In the minimal model program, it is very natural to consider pairs $(X, D)$ where $X$ is a variety and $D$ is a $\bQ$-divisor.  In recent years, the study of pairs $(X, \ba^c)$ where $\ba$ is an ideal sheaf and $c$ is a positive real number, has also become quite common.  Thus it is very natural to try to extend the notion of rational singularities to pairs.  We define two notions of rational pairs.  First we define a \emph{rational pair} which is analogous to a Kawamata log terminal (klt) pair, and then we define a \emph{purely rational triple} which is analogous to purely log terminal (plt) triple (we will discuss the characteristic $p$ analogues later). It is hoped that these definitions and their study will help further the understanding both of rational singularities and log terminal pairs.

In characteristic zero, defining rational singularities for pairs has one distinct advantage over the corresponding variants of log terminal singularities.  In order for $(X, D)$ to be log terminal, one necessarily must have $K_X + D$ a $\bQ$-Cartier divisor.  Likewise, for the pair $(X, \ba^c)$ to be log terminal, $X$ must necessarily be $\Q$-Gorenstein.  One can define rational singularities for a pair $(X, \ba^c)$ without any such conditions on $X$.

Virtually all standard properties of rational singularities transfer to pairs, as we show.  In particular, summands and deformations behave well, see Corollary \ref{rationalSummands} and Theorem \ref{rationalDeforms}, as do various implications between log terminal and rational pairs, see Proposition \ref{rationalGorensteinImpliesklt} and Proposition \ref{kltImpliesrational}.  For the most part, the proofs are generalizations of proofs of the analogous properties of rational singularities.  Since singularities of pairs come up very naturally in theorems related to adjunction and inversion of adjunction, we prove that several of these results extend to rational pairs as well.  In particular, we are able to prove a ``rational'' analogue of inversion of adjunction for log terminal pairs; see Theorem \ref{purelyRationalInversionOfAdjunction}.  Using a similar technique, we are able to give a remarkably short proof of an analogue of inversion of adjunction on log canonicity, which uses the notion of Du Bois singularities; see Theorem \ref{rationalDuBoisInversionOfAdjunction}.

Since the early 1980s, it has been known that rational singularities are closely related to singularities defined by the action of Frobenius map in positive characteristic; see \cite{FedderFPureRat}.  After the introduction of tight closure by Hochster and Hunkeke, see \cite{HochsterHunekeTC1}, a true characteristic $p$ analogue of rational singularities, $F$-rationality, was defined; also see \cite{FedderWatanabe}.  In the next decade, it was shown that a variety has rational singularities if and only if a generic positive characteristic model has $F$-rational singularities;  see \cite{SmithFRatImpliesRat}, \cite{HaraRatImpliesFRat}, and \cite{MehtaSrinivasRatImpliesFRat}.  Thus, we also define $F$-rationality for pairs.  Directly in positive characteristic, we are able to show that $F$-rational pairs satisfy many of the same basic properties that rational pairs do in characteristic zero; see Propositions \ref{basic2}, \ref{ComparisonOfFDivisorialRationalsAndFRegulars}, \ref{inversion} and \ref{SummandTheoremForPositiveCharacteristic} as well as Theorem \ref{FInjectiveInversion}.  Furthermore, building on the techniques of Hara and Yoshida \cite{HaraYoshidaGeneralizationOfTightClosure}, we are able to show a direct correspondence between $F$-rational and rational pairs; see Theorem \ref{correspondence}.

We also relate this to a notion that has existed for many years and was defined and studied for pairs in the toric setting by Manuel Blickle, the \emph{multiplier submodule}; see \cite{SmithTestIdeals}, \cite{HyrySmithNonvanishing}, \cite{BlickleMultiplierModules} and \cite{HaraInterpretation}.  Multiplier ideals and generalized test ideals (their positive characteristic analog) have been studied extensively in recent years as a very powerful invariant which measures singularities of pairs.  For example, a pair is Kawamata log terminal (respectively $F$-regular) if and only if the corresponding multiplier ideal (respectively generalized test ideal) is the entire ring.  When formulating rational singularities associated to pairs, instead of a (multiplier) ideal, it is natural to consider a submodule of the canonical module, an object called the \emph{multiplier submodule} (their characteristic $p$-analogue has been studied under the name ``parameter test submodule''); see definitions \ref{MultiplierSubmoduleDefinition} and \ref{testModuleDefinition}.  Many questions asked about multiplier ideals can also be asked about multiplier submodules; in particular, we look at an analogue of the log canonical threshold in both characteristic zero and positive characteristic; see Definitions \ref{RationalThresholdDefinition} and \ref{DefinitionFInjectiveThreshold}.  We also define jumping exponents for generalized parameter test submodules and show that these numbers form a discrete set of rational numbers under certain conditions; see Definition \ref{DefinitionJumpingExponentForTestModules} and Corollary \ref{rational}.

Most of the techniques in this paper are not new.  They are either techniques related to rational and $F$-rational singularities, or techniques related to log-terminal and $F$-regular singularities.  On the other hand, one might view the fact that these techniques extend so naturally to the cases we consider as further evidence that this generalization of rational singularities to pairs is a natural one.

\begin{rem}
This version of the paper corrects a mistake in the characteristic $p > 0$ theory that appeared in the published version.
In the published version of this paper, we asserted that a version of condition (4) from Lemma \ref{local characterization} was true without the assumption that $\ba$ was principal.  We did not provide a proof and this assertion is not correct.  The problem is that the socle of $H^d_{\bm}(R)$ is not necessarily one-dimensional.  Likewise, Lemma \ref{divisorial characterization}(4) and Remark 6.7 have also been amended.  The proof of Proposition \ref{basic2} has been altered as well (although the new version is no harder).  Finally, this problem also appears with the original definition of $F$-injective pairs we introduced.  However, since we only proved results about $F$-injective pairs $(R, \ba^t)$ when $\ba = (f)$ is a principal ideal, we have restricted the definition of $F$-injective pairs to the case when $\ba$ was principal, see Definition \ref{DefinitionFInjectiveThreshold}.
\end{rem}

{\it Acknowledgements:}
\\
The authors first began discussing the concept of rational singularities of pairs during a workshop held at the American Institute of Mathematics.  The authors would also like to thank the referee for pointing out several typos and in particular thank the referee for pointing out an omission in the definition of purely log terminal / rational triples.

\section{Preliminaries in characteristic zero}
\label{SectionPreliminaries}

All schemes in this paper will be assumed to be separated, noetherian and of essentially finite type over a field.  If $Y$ is a scheme, we will often work in the derived category of $\O_Y$-modules, denoted by $D(Y)$.  The symbol $D^b(Y)$ (respectively $D^+(Y)$, $D^-(Y)$) will denote the derived category of bounded (respectively bounded below, bounded above) complexes of $\O_Y$-modules, $D_{\coherent}(Y)$ (respectively $D_{\quasicoherent}(Y)$) will denote the category of complexes of $\O_Y$-modules with coherent (respectively quasi-coherent) cohomology; see \cite{HartshorneResidues}.  In the setting of the derived category, we will write $F^{\mydot} \qis G^{\mydot}$ if $F^{\mydot}$ and $G^{\mydot}$ are quasi-isomorphic.  We will use $h^i(F^{\mydot})$ to denote the $i$th cohomology of $F^{\mydot}$.  The symbol $\omega_Y^{\mydot}$ will be used to denote a normalized dualizing complex on $Y$, see \cite{HartshorneResidues}, and $\omega_Y$ will be used to denote $h^{-\dim Y}(\omega_Y^{\mydot})$.

We now state \emph{Grothendieck duality} for proper morphisms.

\begin{thm}\cite[III.11.1, VII.3.4]{HartshorneResidues}
\label{GrothendieckDuality}
Let $f : X \rightarrow Y$ be a proper morphism of noetherian schemes of finite dimension.  Suppose $\sF^{\mydot} \in D^{-}_{\quasicoherent}(X)$ and $\sG^{\mydot} \in D^{+}_{\coherent}(Y)$.  Then the duality morphism
\[
\myR f_* \myR \sHom_{X}^{\mydot}(\sF^{\mydot}, f^{!} \sG^{\mydot}) \rightarrow \myR \sHom_{Y}^{\mydot}(\myR f_* \sF^{\mydot}, \sG^{\mydot}),
\]
is an isomorphism.
\end{thm}

\begin{rem}
The case we will consider is when $\sG^{\mydot}$ is a dualizing complex for $Y$ and the map $f$ is a morphism of schemes of finite type over a field $k$ so that $f^{!}(\omega_Y^{\mydot}) = \omega_X^{\mydot}$, giving us the following form of duality
\[
\myR f_* \myR \sHom_{X}^{\mydot}(\sF^{\mydot}, \omega_X^{\mydot}) \cong \myR \sHom_{Y}^{\mydot}(\myR f_* \sF^{\mydot}, \omega_Y^{\mydot}).
\]
\end{rem}

Now we define pairs, log resolutions, and some of the types of characteristic zero singularities we will be considering; see \cite{KollarSingularitiesOfPairs} or \cite{KollarMori} for a more detailed introduction to these definitions.  We fix $X$ to be a noetherian scheme of finite type over a field of characteristic zero $k$.

\begin{defn}
A \emph{pair} $(X, \ba^c)$ is the combined data of a reduced scheme $X$, an ideal sheaf $\ba$ on $X$, and a nonnegative rational (or even real) number $c$.  If $Z$ is a closed subscheme of $X$ defined by an ideal sheaf $I_Z$, then we will often write $(X, cZ)$ instead of the pair $(X, I_Z^c)$.
\end{defn}

\begin{defn}
Suppose that $X$ is as above.  A \emph{resolution} of $X$ is a proper birational map $\pi : \tld X \rightarrow X$ such that $\tld X$ is smooth over $k$.  We let $\exc(\pi)$ denote the exceptional set of $\pi$.  If $\ba$ is an ideal sheaf on $X$, a \emph{log resolution of $\ba$ in $X$} (or simply a \emph{log resolution of $(X, \ba)$} or even a \emph{log resolution of $\ba$}) is a resolution of $X$ such that $\ba \O_{\tld X} = \O_{\tld X}(-G)$ is an invertible sheaf and such that $\exc(\pi) \cup \Supp(G)$ is a simple normal crossings divisor.
\end{defn}

\begin{defn}
A reduced scheme $X$ is said to have \emph{rational singularities} if, for one resolution of $X$, $\pi : \tld X \rightarrow X$, the natural map $\O_X \rightarrow \myR \pi_* \O_{\tld X}$ is a quasi-isomorphism.
\end{defn}

\begin{rem}
If $X$ has rational singularities, then $\O_X \qis \myR \pi_* \O_{\tld X}$ for every resolution.  See for example, \cite[5.10]{KollarMori} or \cite[11.11]{KollarSingularitiesOfPairs}.
\end{rem}

\begin{rem}
It is clear from the definition that rational singularities are necessarily normal.  It also follows immediately from Grauert Riemenschneider vanishing, \cite{GRVanishing}, and Grothendieck duality, see Theorem \ref{GrothendieckDuality}, that rational singularities are Cohen-Macaulay.
\end{rem}

Suppose that $X$ is a normal equidimensional $\bQ$-Gorenstein scheme.  Let $\ba$ be an ideal sheaf on $X$ and suppose that $\pi : \tld X \rightarrow X$ is a log resolution of $(X, \ba^c)$ with $\ba \O_{\tld X} \cong \O_{\tld{X}}(-G)$.  Suppose that $n K_X$ is Cartier, we then define $\pi^*(K_{X})$ to be ${1 \over n}(\pi^* (n K_X))$, which is a $\bQ$-divisor on $\tld X$.  We use $K_{\tld X / X}$ to denote the unique $\bQ$-divisor on $\tld X$, numerically equivalent to $K_{\tld X} - \pi^*(K_X)$ and supported on the exceptional set of $\pi$.

We can now write
\[
K_{\tld X / X} - c G = \sum_{i = 1}^n a(X, E_i, \ba^c) E_i
\]
where the $a(X, E_i, \ba^c)$ are rational numbers and the $E_i$ are divisors.

\begin{defn}
\label{kltDefinition}
The number $a(X, E_i, \ba^c)$ is called the \emph{discrepancy of $(X, \ba^c)$ along the divisor $E_i$}.  We say that $(X, \ba^c)$ has \emph{Kawamata log terminal singularities}, or is simply \emph{klt} if, for a fixed log resolution $\pi$ as above, all of the $a(X, E_i, \ba^c)$ are strictly bigger than $-1$.
\end{defn}

\begin{rem}
The definition of klt singularities is independent of the choice of log resolution; see \cite{KollarMori}.  In fact, if we view each $E_i$ in $\tld X$ as corresponding to a discrete valuation of the fraction field of $X$, then the numbers $a(X, E_i, \ba^c)$ are also independent of the choice of resolution.
\end{rem}

\begin{defn}
\label{multiplierIdealDefinition}
With notation as above, the \emph{multiplier ideal of the pair $(X, \ba^c)$}, denoted by $\mJ(X, \ba^c)$, is defined to be $\pi_* \O_{\tld X}(\lceil K_{\tld X / X} - c G \rceil) \subseteq \O_X$. \end{defn}

\begin{rem}
\label{RemarkMultiplierIdealKLTDefinition}
Note that $(X, \ba^c)$ is klt if and only if $\O_{\tld X}$ is naturally a subsheaf of $\O_{\tld X}(\lceil K_{\tld X / X} - c G \rceil)$.  Thus we see that that $(X, \ba^c)$ is klt if and only if $\mJ(X, \ba^c) = \O_X$
\end{rem}

\begin{rem}
In a context similar to multiplier ideals, we will also often deal with restricting simple normal crossing divisors to a smooth component.  In particular, we will often use the fact that round-down commutes with such restriction without any comment; see \cite[Section 9.1]{LazarsfeldPositivity2}
\end{rem}

\begin{rem}
One can also define log terminal singularities and multiplier ideals for a triple $(X, \Delta, \ba^c)$ where $\Delta$ is a $\bQ$-divisor such that $K_X + \Delta$ is $\bQ$-Cartier.  We will not consider such definitions since this notion does not seem as natural for rational singularities.
\end{rem}

A key property of multiplier ideals that we will rely on is local vanishing, see \cite{EinMultiplierIdealsVanishing}, which is essentially a corollary of Kawamata-Viehweg vanishing; see \cite{KawamataVanishing} and \cite{ViehwegVanishingTheorems}.  We state a formulation of local vanishing for multiplier ideals below.

\begin{thm} [{\cite[9.4]{LazarsfeldPositivity2}}]
\label{MultiplierIdealVanishing}
Using the notation from \ref{multiplierIdealDefinition}, we have
\[
R^j \pi_* \O_{\tld X}(\lceil K_{\tld X /X} - c G \rceil) = 0 \text{  for } j > 0.
\]
\end{thm}

Another variation on log terminal singularities are purely log terminal singularities.  We consider the situation of a triple $(X, H; \ba^c)$ where $X$ is a normal $\bQ$-Gorenstein scheme, $H$ a reduced integral Cartier divisor with ideal sheaf $I_H$, $\ba$ another ideal sheaf and $c$ a nonnegative real number.  A log resolution of such a triple is a simultaneous log resolution of $I_H$ and $\ba$ that is also an embedded resolution of $H$ (which is to say, the strict transform of $H$ is smooth).

\begin{defn}
\label{pltDefinition}
Then we say that such a triple $(X, H;\ba^c)$ has \emph{purely log terminal singularities}, or is simply \emph{plt}, if all the discrepancies of the triple $(X, I_H \ba^c)$ are greater than $-1$, except for those corresponding to the strict transform of $H$ (which are necessarily equal to $-1$).
\end{defn}

\begin{defn}
\label{adjointIdealDefinition}
Let $X$ be a normal $\bQ$-Gorenstein scheme, $H$ a reduced integral Cartier divisor with ideal sheaf $I_H$, $\ba$ another ideal sheaf and $c$ a nonnegative real number.  We define the \emph{adjoint ideal of $(X, H;\ba^c)$}, denoted $\adj(X, H;\ba^c)$ as follows.  Let $\pi : \tld X \rightarrow X$ be a log resolution of $I_H$ and $\ba$ such that the strict transform $\tld H$, of $H$, is smooth (that is, a log resolution of $(X, H; \ba^c)$).  Let $G$ denote the divisor on $\tld X$ such that $\ba \O_{\tld X} = \O_{\tld X}(-G)$.  Then $\adj(X, H;\ba^c)$ is defined to be $\pi_* \O_{\tld X}(\lceil K_{\tld X / X} - c G - \pi^* H + \tld H \rceil) \subseteq \O_X$.
\end{defn}

\begin{rem}
\label{RemarkAdjointIdealPLTDefinition}
Note that $(X, H;\ba^c)$ is plt if and only if $\adj(X, H;\ba^c) = \O_{X}$.
\end{rem}

In the case that $H$ is a Weil-divisor and not a Cartier divisor, one can often still define plt singularities and adjoint ideals for the triple $(X, H;\ba^c)$, (in fact, even further generalizations can be made).  We restrict ourselves to the Cartier case since rational singularities seem best behaved in this context; see remark \ref{RemarkOnNonCartierCase} for additional discussion.

We conclude with a definition of Du Bois singularities; see \cite{DuBoisMain} and \cite{SchwedeEasyCharacterization}.

\begin{defn}
Suppose that $X$ is a reduced scheme embedded as a closed subscheme of a scheme $Y$ with rational singularities.  Let $\pi : \tld Y \rightarrow Y$ be a log resolution of $(Y, X)$ that is an isomorphism outside of $X$ (such log resolutions exist if and only if $Y \backslash X$ is smooth).  Let $E$ denote $(\pi^{-1}(X))_{\red}$.  Then $X$ is said to have \emph{Du Bois singularities} if the natural map $\O_X \rightarrow \myR \pi_* \O_{E}$ is a quasi-isomorphism.
\end{defn}

\begin{rem}
This definition is independent of the choice of embedding or resolution and furthermore, the object $\myR \pi_* \O_{E}$ is also often denoted by $\DuBois{X}$.
\end{rem}

The condition that $\pi$ is an isomorphism outside of $X$ is unnecessary as the following proposition shows; compare with \cite[4.9]{SchwedeEasyCharacterization}.

\begin{prop}
Suppose that $X$ is a reduced closed subscheme of a scheme $Y$ with rational singularities and that $Y \backslash X$ is smooth.  Let $\pi : \tld Y \rightarrow Y$ be a log resolution of the pair $(Y, I_X)$ and let $F$ denote $(\pi^{-1}(X))_{\red}$.  Then $X$ has \emph{Du Bois singularities} if and only if the natural map $\O_X \rightarrow \myR \pi_* \O_{F}$ is a quasi-isomorphism.
\end{prop}

\begin{proof}
It is sufficient to show that $\myR \pi_* \O_{F}$, (or equivalently that $\myR \pi_* \O_{\tld Y}(-F)$) is independent of the choice of resolution.  Since any two log resolutions can be dominated by a third, it is sufficient to consider two log resolutions $\pi_1 : Y_1 \rightarrow Y$ and $\pi_2 : Y_2 \rightarrow Y$ and a map between them $\rho : Y_2 \rightarrow Y_1$ over $Y$.  Let $F_1 = (\pi_1^{-1}(X))_{\red}$ and $F_2 = (\pi_2^{-1}(X))_{\red} = (\rho^{-1}(F_1))_{\red}$.  As mentioned, it is sufficient to prove that $\O_{Y_1}(-F_1) \rightarrow \myR \rho_* \O_{Y_2}(-F_2)$ is a quasi-isomorphism.  Dualizing the map and applying Grothendieck duality implies that it is sufficient to prove that $\omega_{Y_1}(F_1) \leftarrow \myR \rho_* ( \omega_{Y_2}(F_2))$ is a quasi-isomorphism.

We now apply the projection formula while twisting by $\omega_{Y_1}^{-1}(-F_1)$ (which is invertible since $Y_1$ is smooth).  Thus it is sufficient to prove that
\[
\myR \rho_* (\omega_{Y_2/Y_1}(F_2 - \rho^* F_1)) \rightarrow \O_{Y_1}
\]
is a quasi-isomorphism.  But note that $F_2 - \rho^* F_1 = - \lfloor \rho^* (1-\epsilon) F_1 \rfloor$ for sufficiently small $\epsilon > 0$.  Thus it is sufficient to prove that the pair $(Y_1, (1-\epsilon)F_1)$ has klt singularities by local vanishing for multiplier ideals; see \cite[9.4]{LazarsfeldPositivity2}.  But this is true since $Y_1$ is smooth and $F_1$ is a reduced integral divisor with simple normal crossings.  Compare this proof with the proof of Theorem \ref{rationalDuBoisInversionOfAdjunction}.
\end{proof}

\begin{rem}
While it is hoped that the condition that $Y \backslash X$ is smooth can be removed, see \cite{SchwedeEasyCharacterization}, it follows from \cite{KovacsDuBoisLC1} that if $\O_X \rightarrow \myR \pi_* \O_{F}$ is a quasi-isomorphism (for any $Y$, even without rational singularities), then $X$ has Du Bois singularities.
\end{rem}

\section{Basic definitions and fundamental properties in characteristic zero}
\label{SectionDefinitionInCharacteristicZero}

\begin{defn}
\label{rationalDefinition}
Let $(X, \ba^c)$ be a pair and let $\pi : \tld X \rightarrow X$ with $\ba \O_{\tld X} = \O_{\tld X}(-G)$ be a log resolution of $\ba$.  We say that the pair $(X, \ba^c)$ has \emph{rational singularities} (or \emph{Kawamata rational singularities}) if the natural map $\O_X \rightarrow \myR \pi_* \O_{\tld X}(\lfloor c G \rfloor)$ is a quasi-isomorphism.
\end{defn}

\begin{rem}
Explicitly, the pair $(X, \ba^c)$ has rational singularities if and only if $\O_X \rightarrow \pi_* \O_{\tld X}(\lfloor c G \rfloor)$ is an isomorphism and $R^i \pi_* \O_{\tld X}(\lfloor c G \rfloor) = 0$ for $i > 0$.  Also note that the natural map of $\O_X$ to its normalization can composed with the map $\pi_* \O_{\tld X} \rightarrow \pi_* \O_{\tld X}(\lfloor c G \rfloor)$ to obtain $\O_X \rightarrow \pi_* \O_{\tld X}(\lfloor c G \rfloor)$, proving that $\O_X$ is a summand of its normalization and is thus normal.
\end{rem}

\begin{rem}
\label{DualrationalFormulation}
By Grothendieck duality, the pair $(X, \ba^c)$ has rational singularities if and only if the natural map $\myR \pi_* \omega_{\tld X}^{\mydot} \tensor \O_{\tld X}(\lceil - c G \rceil) \rightarrow \omega_X^{\mydot}$ is an isomorphism.  Compare with \cite[Page 50]{KempfToroidalEmbeddings}.
\end{rem}

Our first goal is to prove that this definition is independent of the choice of resolution.

\begin{prop}
\label{RationalInPairsIndependent}
The definition given in \ref{rationalDefinition} is independent of the choice of resolution.
\end{prop}
\begin{proof}
Let $(X, \ba^c)$ be a pair as in \ref{rationalDefinition}.  Since any two log resolutions can be dominated by a third, it is enough to consider two log resolutions of $\ba$, $X_1$ and $X_2$ with a map between them, as pictured below.
\[
\xymatrix{
X_2 \ar[rr]^{\rho} \ar[dr]_{\pi_2} & & X_1 \ar[dl]^{\pi_1} \\
& X & \\
}
\]
Let us use the $G_1$ and $G_2$ to denote the divisors (on $X_1$ and $X_2$ respectively) such that  $\ba \O_{X_1} = \O_{X_1}(-G_1)$ and $\ba \O_{X_2} = \O_{X_2}(-G_2)$.  It is enough to prove that the map $\O_{X_1}(\lfloor c G_1 \rfloor) \rightarrow \myR \rho_* \O_{X_2}(\lfloor c G_2 \rfloor)$ is a quasi-isomorphism (note such a map exists since $\rho^* \lfloor cG_1 \rfloor \leq \lfloor cG_2 \rfloor$).  By Grothendieck duality (since $X_1$ and $X_2$ are smooth), this is equivalent to proving that we have a quasi-isomorphism
\[
\myR \rho_* \omega_{X_2}(-\lfloor c G_2 \rfloor) \rightarrow \omega_{X_1}(-\lfloor c G_1 \rfloor).
\]
Tensoring this map with $\tensor \omega_{X_1}^{-1}$ (which is an invertible sheaf since $X_1$ is smooth), then reduces our question to independence of the definition of multiplier ideals (after an application of local vanishing for multiplier ideals, \cite[9.4]{LazarsfeldPositivity2}), since $\rho^* G_1 = G_2$ and $G_1$ is a simple normal crossings divisor.  Also see \cite[Secition 2]{GRVanishing}.
\end{proof}

Our next main goal is to explore how varying the constant $c$ or varying the ideal $\ba$ affects whether the pair in question has rational singularities.  In the process of doing this, we will introduce a notion analogous to the multiplier ideal and will also prove a technical theorem, \ref{GeneralizedKovacsSplittingTheorem}, related to \cite[Theorem 1]{KovacsRat} and \cite[5.13]{KollarMori}, which will be used to give a simple proof that log terminal pairs are rational and that summands of (appropriate) rational pairs are rational.
The essential ingredient in all of this is the following (vanishing) lemma.  This lemma, which will be obvious to experts, can be thought of as either a generalization of Grauert Riemenschneider vanishing, see \cite{GRVanishing}, or a slight modification of the usual formulation of local vanishing for multiplier ideals, \cite[9.4.1]{LazarsfeldPositivity2}.

\begin{lem}
\label{LocalVanishingForMultiplierSubmodules}
Suppose that $X$ is a reduced equidimensional scheme and $\ba$ is an ideal sheaf on $X$.  Further suppose that $\pi : \tld X \rightarrow X$ is a log resolution of $\ba$ with $\ba \O_{\tld X} = \O_{\tld X}(-G)$.  Then for any nonnegative real number $c$ and for all $i > 0$ we have $h^i(\myR (\pi_* \omega_{\tld X} \tensor \O_{\tld X}(\lceil - c G \rceil)))  = 0$.
\end{lem}
\begin{proof}
First note that we may assume that $X$ is normal since the map $\pi$ factors through the normalization of $X$ and finite maps have no higher cohomology.  Thus, we may also assume that $X$ is irreducible.  We then reduce to the case when $\ba$ is a (locally) principal ideal sheaf by choosing general elements of $\ba$; see \cite[9.2.22-9.2.28]{LazarsfeldPositivity2}.  The proof is then the same as the proof of \cite[9.4.1, 9.4.17]{LazarsfeldPositivity2}, except we do not need to pull back $K_X$.  The essential ingredient is the Kawamata-Viehweg vanishing theorem; see \cite{KawamataVanishing}, \cite{ViehwegVanishingTheorems}.
\end{proof}

In her thesis and in \cite{SmithTestIdeals}, Karen Smith noted that when dealing with rational singularities and related notions, instead of working with (analogues of) multiplier ideals, one should work with submodules of the canonical module.  This idea was further studied in \cite{HaraInterpretation}.  Thus the following definition is natural, also compare with Remark \ref{testModuleDefinition}

\begin{defn}\cite{BlickleMultiplierModules}
\label{MultiplierSubmoduleDefinition}
The \emph{multiplier submodule} of a pair $(X, \ba^c)$ is defined to be the image of $\pi_* (\omega_{\tld X} \tensor \O_{\tld X}(\lceil - c G \rceil))$ inside $\omega_X$.  We will denote it by $\mJ(\omega_X, \ba^c)$.
\end{defn}

It is easy to see that this submodule is independent of the choice of resolution.  From this point of view, lemma \ref{LocalVanishingForMultiplierSubmodules} can be thought of as \emph{local vanishing for multiplier submodules}.

\begin{lem}
\label{InjectivityOfMultiplierSubmoduleLemma}
If $X$ is a reduced equidimensional scheme as above, the natural map $\pi_* (\omega_{\tld X} \tensor \O_{\tld X}(\lceil - c G \rceil)) \rightarrow \omega_X$ is injective.
\end{lem}
\begin{proof}
Consider the exact triangle
\[ \xymatrix{ \O_X \ar[r] & \myR \pi_* \O_{\tld X} \ar[r] & C^{\mydot} \ar[r]^-{+1} &  } \]
and note that $\dim(\Supp(h^i(C^\mydot))) < \dim X - i$.  By an easy analysis of a spectral sequence, one has $\bH^d_{\bm}(C^{\mydot}) = 0$ which implies that there is a surjection $H^d_{\bm}(\O_X) \rightarrow \bH^d_{\bm}(\myR \pi_* \O_{\tld X})$ (for any maximal ideal $\bm$).  By local duality, \cite[V, Theorem 6.2]{HartshorneResidues}, the natural map $\pi_* \omega_{\tld X} \rightarrow \omega_X$ is an injection.  But then we are done since
$\omega_{\tld X}(\lceil - c G \rceil) \subset \omega_{\tld X}$ and $\pi_*$ is left exact.  Compare with \cite[Section 2, Remark (b)]{LipmanTeissierPseudoRational} and with \cite[Section 1]{KempfSomeQuotientVarietiesHaveRationalsSingularities}.
\end{proof}

\begin{cor}
\label{MultiplierSubmoduleCharacterization}
Suppose that $X$ is a reduced equidimensional Cohen-Macaulay scheme and $\ba$ is an ideal sheaf on $X$.  Then $(X, \ba^c)$ has rational singularities if and only if the multiplier submodule of $X$, $\pi_* (\omega_{\tld X} \tensor \O_{\tld X}(\lceil - c G \rceil))$, is equal to $\omega_X$.
\end{cor}

At this point, it is natural to mention a (characteristic-free) definition for rational singularities of pairs that makes sense even when $X$ isn't known to have a resolution.  This slight generalization of a definition of Lipman and Teissier will appear later in the paper when comparing rational and $F$-rational pairs, see Theorem \ref{correspondence}.

\begin{defn}[\textup{cf. \cite[Sec. 2]{LipmanTeissierPseudoRational}}]\label{pseudo-rational}
Let $(R,\m)$ be a $d$-dimensional reduced local ring, let $\a \subseteq R$ be an ideal such that $\a$ contains elements not contained in any minimal prime of $R$ and let $t \ge 0$ be a real number.
Then $(R, \a^t)$ is \textit{pseudo-rational} if $R$ is normal, Cohen-Macaulay, analytically unramified, and if for any proper birational morphism $\pi:Y \to X:=\Spec R$ from a normal scheme $Y$, such that $\a \O_Y=\O_Y(-G)$ is invertible, the map
$$\delta_{\pi}:H^d_{\m}(R) \to H^d_E(\O_Y(\lfloor tG \rfloor))$$
is injective, where $E=\pi^{-1}(\m)$ denotes the closed fiber of $\pi$ and $\delta_{\pi}$ is the map induced by $\O_{\Spec R} \rightarrow \myR \pi_* \O_Y(\lfloor tG \rfloor)$.
\end{defn}

\begin{rem}\label{pseudo rem}
In addition, when $R$ is essentially of finite type over a field of characteristic zero, a straightforward application of local duality (see \cite[V, Theorem 6.2]{HartshorneResidues}) implies that $(R,\a^t)$ is pseudo-rational if and only if $(\Spec R,\a^t)$ has rational singularities.
\end{rem}

Now we come to the promised generalization of a result of Kov\'acs, \cite{KovacsRat}.

\begin{thm}
\label{GeneralizedKovacsSplittingTheorem}
Suppose that $(X, \ba^c)$ is a pair such that $\pi : \tld X \rightarrow X$ is a log resolution of $\ba$.  If the natural map
\[
\O_X \rightarrow \myR \pi_* \O_{\tld X}(\lfloor c G \rfloor)
\]
has a left inverse (meaning that there exists a map $\O_{\tld X}(\lfloor c G \rfloor) \rightarrow \O_X$ such that the composition
$
\O_X \rightarrow \myR \pi_* \O_{\tld X}(\lfloor c G \rfloor) \rightarrow \O_X
$
is a quasi-isomorphism), then $(X, \ba^c)$ has rational singularities.
\end{thm}
The proof is virtually the same as the one found in \cite{KovacsRat}, we simply use \ref{LocalVanishingForMultiplierSubmodules} instead of Grauert Riemenschneider vanishing.
\begin{proof} \cite{KovacsRat}
We first note that since $\pi$ factors through the normalization of $X$, we immediately see that $\O_X$ is a summand of its own normalization, and thus itself normal.  Therefore, without loss of generality, we may assume that $X$ is irreducible (and in particular, equidimensional).
Now, apply Grothendieck duality to give us the composition
\[
\omega_X^\mydot \rightarrow \myR \pi_* \omega_{\tld X}^{\mydot}(-\lfloor c G \rfloor) \rightarrow \omega_X^{\mydot}.
\]
By lemma \ref{LocalVanishingForMultiplierSubmodules}, and since the composition is an isomorphism, we have $h^i(\omega_X^{\mydot}) = 0$ for $i \neq -\dim X$.  This implies that $X$ is Cohen-Macaulay.  It is now enough to show that
\[
\pi_* (\omega_{\tld X} \tensor \O_{\tld X}(-\lfloor c G \rfloor) ) \rightarrow \omega_X
\]
is an isomorphism.  However, the map is injective by \ref{InjectivityOfMultiplierSubmoduleLemma}.  It is surjective since it is a split surjection (by assumption).
\end{proof}

One could have given an indirect argument that $X$ is Cohen-Macaulay by first showing that $X$ is rational, but Kov\'acs' argument generalizes quite well to pairs and is really no longer than an indirect argument.

\begin{cor}
Suppose that $X$ is a reduced scheme, $\ba$ is an ideal sheaf and $c_1 < c_2$ are nonnegative real numbers.  If $(X, \ba^{c_2})$ has rational singularities, so does $(X, \ba^{c_1})$.  Furthermore if $\bb$ is another ideal sheaf with $\ba \subseteq \bb$, then if $(X, \ba^{c_1})$ is rational, so is $(X, \bb^{c_1})$.
\end{cor}
\begin{proof}
Let $\pi : \tld X \rightarrow X$ be a log resolution of $\ba$.  We have the following composition
\[
\O_X \rightarrow \myR \pi_* \O_{\tld X}(\lfloor c_1 G \rfloor) \rightarrow \myR \pi_* \O_{\tld X}(\lfloor c_2 G \rfloor)
\]
which is a quasi-isomorphism by assumption, proving that $(X, \ba^{c_1})$ has rational singularities by \ref{GeneralizedKovacsSplittingTheorem}.  The proof of the second statement is similar.
\end{proof}

\begin{cor}
\label{rationalImpliesRational}
Suppose that the pair $(X, \ba^c)$ has rational singularities, then $X$ has rational singularities and is in particular Cohen-Macaulay.
\end{cor}

\begin{rem}
In the previous two corollaries, one can avoid working with the derived category by first dualizing and considering containments of multiplier submodules.
\end{rem}

We conclude this section with a definition of purely rational singularities, compare with Definition \ref{DivisoriallyFRationalDefinition}.

\begin{defn}
\label{purelyRationalDefinition}
Let $X$ be a normal scheme, $H$ an integral reduced Cartier divisor with ideal sheaf $I_H$, $\ba$ another ideal sheaf with no minimal prime among the components of $H$, and $c$ a nonnegative real number.  Suppose that $\pi : \tld X \rightarrow X$ is a log resolution of $H$ and $\ba$, where $\tld H$, the strict transform of $H$, is smooth (that is, $\pi$ is a log resolution of the triple $(X, H; \ba^c)$).  Let us use $G$ to denote the divisor on $\tld X$ such that $\ba \O_{\tld X} = \O_{\tld X} (-G)$.  Then we say that $(X, H; \ba^c)$ has \emph{purely rational singularities} if the natural map
\[
\O_X \rightarrow \myR \pi_* \O_{\tld X}(\lfloor c G + \pi^* H - \tld H \rfloor)
\]
is a quasi-isomorphism.
\end{defn}

\begin{rem}
By Grothendieck duality, $(X, H; \ba^c)$ has purely rational singularities if and only if
\[
\myR \pi_* (\omega_{\tld X}^{\mydot} \tensor \O_{\tld X}(\lceil-cG - \pi^* H + \tld H \rceil) ) \rightarrow \omega_X^{\mydot}
\]
is a quasi-isomorphism.
\end{rem}

We also define the adjoint submodule.

\begin{defn}
\label{adjointSubmoduleDefinition}
Suppose $X$ is a reduced scheme, $H$ a Cartier divisor and $\ba$ an ideal sheaf with no minimal primes in common with any of the components of $H$.  The \emph{adjoint submodule} of a triple $(X, H;\ba^c)$ is defined to be the image of $\pi_* (\O_{\tld X}(\lceil K_X -cG - \pi^* H + \tld H \rceil)$ inside $\omega_X$ where $\pi$ is defined as above.  We will denote the adjoint submodule by $\adj(\omega_X, H; \ba^c)$.
\end{defn}

We now show that the notions of purely rational singularities and the adjoint submodule are well defined.

\begin{prop}
\label{PurelyRationalWellDefinedAndVanishing}
With the notation as in definition \ref{purelyRationalDefinition}, the definition of purely rational singularities is independent of the choice of resolution (more generally, the \emph{adjoint submodule} $\pi_* (\O_{\tld X}(\lceil K_X -cG - \pi^* H + \tld H \rceil) \subset \omega_X$ is well defined).  Furthermore,
\[
h^i(\myR \pi_* (\O_{\tld X}(\lceil K_{\tld X} -cG - \pi^* H + \tld H \rceil) ) ) = 0
\]
for $i > 0$ so that $(X, H;\ba^c)$ has purely rational singularities if and only if $X$ is Cohen-Macaulay and
\[
\pi_* (\O_{\tld X}(\lceil K_X -cG - \pi^* H + \tld H \rceil)) \rightarrow \O_{\tld X}(K_X)
\]
is surjective (in other words, if and only if the adjoint submodule is equal to $\omega_X$).
\end{prop}
\begin{proof}
To show that $h^i(\myR \pi_* \omega_{\tld X}(-\lfloor c G + \pi^* H - \tld H \rfloor)) = 0$ for $i > 0$, it is enough to show that $h^i(\myR \pi_* \omega_{\tld X}(-\lfloor c G - \tld H \rfloor)) = 0$ by the projection formula since $H$ is Cartier.  Thus consider the following short exact sequence,
\[
0 \rightarrow \omega_{\tld X}(-\lfloor cG \rfloor) \rightarrow \omega_{\tld X}(-\lfloor c G - \tld H \rfloor) \rightarrow \omega_{\tld H}(-\lfloor c G \rfloor) \rightarrow 0.
\]
If we apply $\myR \pi_*$ we see that the higher cohomology of $\myR \pi_* \omega_{\tld X}(-\lfloor cG \rfloor)$ is zero by \ref{LocalVanishingForMultiplierSubmodules} and likewise the higher cohomology of $\myR \pi_* \omega_{\tld H}(-\lfloor c G \rfloor)$ is also zero.  This proves the vanishing we desired.

It is now sufficient to prove that the adjoint submodule is well defined.  Since any two log resolutions can be dominated by a third, we consider the case of two log resolutions of $(X, H; \ba^c)$, $X_1$ and $X_2$ with a map between them, as pictured below.
\[
\xymatrix{
X_2 \ar[rr]^{\rho} \ar[dr]_{\pi_2} & & X_1 \ar[dl]^{\pi_1} \\
& X & \\
}
\]
As before, we assume that the strict transforms, $\tld H_1$ and $\tld H_2$, of $H$ in $X_1$ and $X_2$ are smooth.  We define $G_i$ to be the divisor on $X_i$ such that $\ba \O_{X_i} = \O_{X_i} (-G_i)$ and note that $\rho^* G_1 = G_2$ and $\rho^* \pi_1^* H = \pi_2^* H$.

It is now sufficient to show that the map
\[
\rho_* \omega_{X_2}(-\lfloor c G_2 + \pi_2^* H - \tld H_2 \rfloor) \rightarrow \omega_{X_1}
\]
has image $\omega_{X_1}(-\lfloor c G_1 + \pi_1^* H - \tld H_1 \rfloor)$.  By the projection formula (twisting by $\tensor (\omega_{X_1}^{-1}(\pi_1^* H - \tld H_1))$), it is sufficient to show that the image of the map
\[
\rho_* \O_{X_2} (\lceil K_{X_2/X_1} - c G_2 + \tld H_2 - \rho^* \tld H_1 \rceil) \rightarrow \O_{X_1}(\pi_1^* H - \tld H_1)
\]
is equal to $\O_{X_1}(-\lfloor c G_1\rfloor) = \O_{X_1}(-\lfloor c G_1\rfloor - \tld H_1 + \tld H_1  )$.  But this follows since the definition of the usual adjoint ideal is independent of choice of resolution.
\end{proof}

\begin{rem}
\label{adjointMultiplierShortExactSequence}
With the notation as above, the previous proof implies we have a short exact sequence
\[
0 \rightarrow \mJ(\omega_X, \ba^c) \rightarrow \adj(\omega_X, H; \ba^c) \tensor \O_X(H) \rightarrow \mJ(\omega_H, (\ba|_H)^c) \rightarrow 0.
\]
Note this is essentially the same as \cite[9.3.44]{LazarsfeldPositivity2}.  Also see Theorem \ref{purelyRationalInversionOfAdjunction}.
\end{rem}

Using the vanishing in Proposition \ref{PurelyRationalWellDefinedAndVanishing}, one can prove an analogue of \ref{GeneralizedKovacsSplittingTheorem}.

\begin{thm}
\label{PurelyRationalSplitting}
If the natural map $\O_X \rightarrow \myR \pi_* \O_{\tld X}(\lfloor c G + \pi^* H - \tld H \rfloor)$ has a left inverse (that is, there exists a map $\myR \pi_* \O_{\tld X}(\lfloor c G + \pi^* H - \tld H \rfloor) \rightarrow \O_X$ such that the composition,
\[
\O_X \rightarrow \myR \pi_* \O_{\tld X}(\lfloor c G + \pi^* H - \tld H \rfloor) \rightarrow \O_X,
\]
is a quasi-isomorphism) then $(X, H;\ba^c)$ has purely rational singularities.
\end{thm}

The proof is the same as in \ref{GeneralizedKovacsSplittingTheorem}.

\begin{rem}
\label{PurelyRationalImpliesRationalRemark}
Note that if $(X, H;\ba^c)$ has purely rational singularities, then $(X, \ba^c I_H^{(1-\epsilon)})$ has (Kawamata) rational singularities for every epsilon satisfying $1 \geq \epsilon \geq 0$ by \ref{GeneralizedKovacsSplittingTheorem}; in particular, $(X, \ba^c)$ has rational singularities.
\end{rem}

\begin{rem}
\label{RemarkOnNonCartierCase}
Finally, let us briefly discuss the case where $H$ is not Cartier.  In such a case, one can consider $\pi_* (\O_{\tld X}(\lceil K_X -cG + \tld H \rceil) \subset \omega_X(H)$ instead of the adjoint submodule.  One still has a vanishing for the higher cohomology in such a case, and many results still work.  This object seems somewhat contrived however, and doesn't seem as closely related to the adjoint ideals as defined, for example, in \cite{LazarsfeldPositivity2}.  For this reason, we strict ourselves to the Cartier case.
\end{rem}

\section{Log terminal singularities, deformations, summands, and adjunction}
\label{SectionPropertiesOfRationalPairs}

In this section, we show how rational pairs relate to log terminal pairs, prove that pairs with rational singularities behave well with respect to deformation and summands, and conclude this section by showing that rational pairs satisfy several inversion of adjunction type results often observed for log terminal pairs.  We also give a simple proof of a result related to inversion of adjunction on log canonicity, which uses the notion of Du Bois singularities.

First we relate log terminal and rational singularities associated to pairs.  In particular, we show that Kawamata log terminal pairs are rational and that rational pairs $(X, \ba^c)$ with $X$ Gorenstein, are Kawamata log terminal; also see \cite{ElkikRationalityOfCanonicalSings}.  We then compute an example which shows that these notions are distinct even when $X$ is $\bQ$-Gorenstein.  Compare the following two results with Propositions \ref{basic2} and \ref{ComparisonOfFDivisorialRationalsAndFRegulars}.

\begin{prop}
\label{rationalGorensteinImpliesklt}
Suppose $(X, \ba^c)$ is rational (respectively $(X, H;\ba^c)$ is purely rational) and $X$ is Gorenstein, then $(X, \ba^c)$ is klt (respectively $(X, H;\ba^c)$ is plt).
\end{prop}
\begin{proof}
Let $\pi : \tld X \rightarrow X$ be a log resolution of $\ba$.  By \ref{DualrationalFormulation}, we have a quasi-isomorphism $\myR (\pi_* \omega_{\tld X}^{\mydot} \tensor  \O_{\tld X}(\lceil - c G \rceil)) \qis \omega_X^{\mydot}$.  But then since $\omega_X$ is a line bundle, we have $\pi_*  \O_{\tld X}(\lceil K_{\tld X / X} - c G \rceil) \cong \O_X$ by the projection formula.  Thus the pair is klt.  The proof of the plt case is the same.
\end{proof}

\begin{prop}
\label{kltImpliesrational}
Suppose that $X$ is $\bQ$-Gorenstein and that $(X, \ba^c)$ is klt (respectively $(X, H;\ba^c)$ is plt), then $(X, \ba^c)$ is also rational (respectively $(X, H;\ba^c)$ is purely rational).
\end{prop}
\begin{proof}
This statement is local, so we may assume that $X$ is affine.  Let $\pi : \tld X \rightarrow X$ be a log resolution of $\ba$.  Now, since $(X, \ba^c)$ is klt, we have a natural inclusion
\[
\O_{\tld X} \subseteq \O_{\tld X}(\lceil K_{\tld X / X} - c G \rceil).
\]
This implies an inclusion
\[
\O_{\tld X}(\lfloor cG \rfloor) \subseteq \O_{\tld X}(\lceil K_{\tld X / X} - c G \rceil + \lfloor cG \rfloor) \subseteq \O_{\tld X}(\lceil K_{\tld X / X} \rceil).
\]
Applying $\myR \pi_*$ gives us a composition
\[
\O_X \rightarrow \myR \pi_* \O_X(\lfloor cG \rfloor) \rightarrow \myR \pi_* \O_{\tld X}(\lceil K_{\tld X / X} \rceil).
\]
But $\myR \pi_* \O_{\tld X}(\lceil K_{\tld X / X}\rceil)$ is quasi-isomorphic to $\O_X$ since $X$ is log terminal (using local vanishing for multiplier ideals, \cite[9.4]{LazarsfeldPositivity2}), which completes the proof of the klt case by \ref{GeneralizedKovacsSplittingTheorem}.

In the plt case, the proof is analogous.  We begin with the inclusion $\O_{\tld X} \subseteq \O_{\tld X}(\lceil K_{\tld X/X} - cG - \pi^* H + \tld H\rceil)$ and observe that $\pi^* H - \tld H$ is a integral divisor.  This gives us an inclusion $\O_{\tld X}(\lfloor cG + \pi^*H - \tld H \rfloor) \subseteq \O_{\tld X}(\lceil K_{\tld X / X} \rceil)$ where $\tld H$ is the strict transform of $H$.  We then apply $\myR \pi_*$ and use \ref{PurelyRationalSplitting} which completes the proof.
\end{proof}

\begin{rem}
Compare the previous proof with \cite[Theorem 4]{KovacsRat}.
\end{rem}

\begin{rem}
If $X$ is not $\bQ$-Gorenstein, but $((X, \Delta), \ba^c)$ is klt (in particular, $K_X + \Delta$ is $\bQ$-Cartier), then the same proof implies that $(X, \ba^c)$ is rational.  In such a case, it would seem natural to try to show that $(X, \Delta)$ is a rational pair, however, there is no clear way to pull back $\Delta$ as a divisor, since it is not $\bQ$-Cartier by assumption.
\end{rem}

\begin{rem}
One could also give a more indirect proof of Theorem \ref{kltImpliesrational} which is less homological by comparing multiplier ideals and multiplier submodules.
\end{rem}

We now present an example of a pair with a log terminal underlying scheme, which has rational and not log terminal singularities.

\begin{eg}
\label{RationalNotKLTExample}
Consider the surface singularity $X = \Spec \bC[x^3, x^2 y, x y^2, y^3]$.  This is a surface with cyclic quotient singularities, and so it is in particular, log terminal. First let us consider this scheme's resolution and how this affects its canonical divisor.  This singularity can be resolved by a single blow-up $\pi : \tld X \rightarrow X$ at the ideal $\bm = (x^3, x^2 y, x y^2, y^3)$.

The canonical module $\omega_X$ can be identified with the ideal $(x^2 y, x y^2)$.  Note that with this identification, we have $\omega_X^{(3)} \cong (x^3 y^3)$.  Now, it is easy to see that the pair $(X, \bm^t)$ is rational for any $0 < t < 1$, but not klt for $t$ sufficiently close to $1$ since $X$'s discrepancy along the exceptional divisor is equal to $-{1 \over 3}$.
\end{eg}

The previous example also suggests the following definition.  As an analogy with the log canonical threshold, one can define the follow rational number, compare with Definition \ref{DefinitionFInjectiveThreshold}.

\begin{defn}
\label{RationalThresholdDefinition}
Let $X$ be a scheme with rational singularities and $\ba$ an ideal sheaf.  We define the \emph{rational threshold} of the pair $(X, \ba)$, denoted by $\text{rt}(X, \ba)$, to be equal the following number:
\[
\text{rt}(X, \ba) = \sup\{t > 0 | (X, \ba^t) \text{ has rational singularities} \}
\]
\end{defn}

In the Example \ref{RationalNotKLTExample}, the log canonical threshold of the pair $(X, \ba)$ is equal to ${2 \over 3}$, whereas the rational threshold is equal to $1$.  More generally, suppose that $X$ is a variety with a log resolution $\pi : \tld X \rightarrow X$ which has only a single reduced exceptional divisor $E$ which dominates $P$ and was obtained by blowing up the same ideal $P$ where $P \O_{\tld X} = \O_{\tld X}(-E)$ , then the rational threshold of $(X, P)$ is always an integer.  On the other hand, there are many examples of varieties with non-integer rational thresholds since the rational threshold and the log canonical threshold of a Gorenstein scheme clearly coincide.

Let us consider now a broader set of examples, that of the Veronese subrings.  We will use a slightly different approach from the example above.  The following generalization of a lemma by Kov\'acs will be useful in this computation.

\begin{lem} \cite[Lemma 3.3]{KovacsDuBoisLC1} \label{LemmaMuchCohomologyVanishes}
Suppose that $X$ is a Cohen-Macaulay scheme of essentially finite type over a field of characteristic zero.  Suppose that $\ba$ is an ideal sheaf and $t$ is a positive rational number.  Let $\Sigma$ be the subset of $X$ where $(X, \ba^t)$ does not have rational singularities.   Let $\pi : \tld X \rightarrow X$ is a log resolution of $(X, \ba^t)$ with $\ba \O_{\tld X} = \O_{\tld X}(-G)$.  Then $R^i (\O_{\tld X}(\lfloor tG \rfloor)) = 0$ for all $0 < i < \dim X - \dim \Sigma - 1$.
\end{lem}
\begin{proof}
The proof is virtually the same as in \cite{KovacsDuBoisLC1}, one simply uses Theorem \ref{LocalVanishingForMultiplierSubmodules} instead of Grauert Riemenschneider vanishing.
\end{proof}

\begin{eg}
Suppose $S = k[x_1, \ldots, x_d]$.  Let $R$ be the $r$th Vernonese subring, $R = k[x_1^r, x_1^{r-1} x_2, x_1^{r-1}x_2^2, \ldots, x_{d-1} x_d^{r-1}, x_d^r]$.  We are going to study the rational threshold of the pair $(\Spec R, \bm^t)$ where $\bm$ is the maximal ideal of the origin.  It is clear that the pair can be resolved with a single blow-up, and to study that blow-up, we can use a set of $d$ charts corresponding to placing each $x_i^r$ in the denominator.  Note that this implies that the rational threshold must be an integer.  Fix $X = \Spec R$ and $\pi : \tld X \rightarrow X$ the aforementioned resolution and let $E$ be the exceptional divisor (note $\bm \O_{\tld X} = \O_{\tld X}(-E)$).  Since $R$ is a Cohen-Macaulay isolated singularity, it is sufficient to understand the cohomology $R^{d-1} \pi_* \O_{\tld X}(\lfloor t E \rfloor)$ by Lemma \ref{LemmaMuchCohomologyVanishes}.

We use \Cech cohomology to interpret this object.  Using the charts corresponding to $x_1^r, \ldots, x_d^r$, we see that an arbitrary element of $R^{d-1} \pi_* \O_{\tld X}(\lfloor t E \rfloor)$ looks like $f \over {(x_1^r x_2^r \ldots x_d^r)^c}$.  Note that the order of vanishing of $f$ along $E$ must be greater than or equal to $c d - \lfloor t \rfloor$.  A natural first nonzero element of the cohomology group would seem to be $ {(x_1^{r-1} x_2^{r-1} \ldots x_d^{r-1})} \over {(x_1^r x_2^r \ldots x_d^r)}$, unfortunately, that numerator doesn't always exist in this context since in some sense the numerator's order of vanishing on  $E$ is $d(r-1)/r$, which is not always an integer.  In order to see that the pair is non-rational it is natural to seek a cohomology element which vanishes on $E$ to degree $\lfloor d(r-1)/r \rfloor - d = \lfloor -d/r \rfloor$, which by assumption must be greater than or equal to $-\lfloor t \rfloor$.  It is not hard to see that such a non-zero element exists assuming the arithmetic is satisfied by modifying the original ``non-existent'' element above.  In other words, if $t \geq \lceil d/r \rceil$, then $(X, \bm^t)$ cannot be a rational pair, which means that $\text{rt}(X, \bm) \leq \lceil d / r \rceil$.  On the other hand, the log canonical threshold $\text{lt}(X, \bm)$ is equal to $d / r$, and it is clear that $\text{lt}(X, \bm) < \text{rt}(X, \bm)$.  Thus we have the following inequality,
$d/r \leq \text{rt}(X, \bm) \leq \lceil d / r \rceil$.
Therefore, since $\text{rt}(X, \bm)$ is an integer, it must be equal to $\lceil d / r \rceil$.
\end{eg}

See Example \ref{VeroneseExamples}, for a study of the same class of singularities using positive characteristic techniques (an explicit proof of the fact that $\text{rt}(X, \bm) \geq \lceil d / r \rceil$ is given in positive characteristic).

We now prove that summands of pairs with rational singularities are rational; compare with \cite{BoutotRational}.  In fact, we prove a more general result analogous to the full generality of \cite[Theorem 1]{KovacsRat}.  The proof is relatively short, the key ingredient is Theorem \ref{GeneralizedKovacsSplittingTheorem}, and was inspired by a similar result of \cite{KovacsRat}.

\begin{thm}
\label{GeneralizedKovacsMorphismSplitting}
Suppose that $\rho : Y \rightarrow X$ is a dominant morphism of reduced schemes such that every component of $Y$ dominates a component of $X$.  Let $\ba$ be an ideal sheaf on $X$ and suppose that $(Y, (\ba \O_Y)^c)$ is rational.  Further suppose that the natural map $\O_X \rightarrow \myR \rho_* \O_Y$ has a left inverse (that is, there exists a map $\delta : \myR \rho_* \O_Y \rightarrow \O_X$ such that $\O_X \rightarrow \myR \rho_* \O_Y \rightarrow \O_X$ is a quasi-isomorphism).  Then $(X, \ba^c)$ has rational singularities as well.
\end{thm}
\begin{proof}
Let $\pi : \tld X \rightarrow X$ and $\pi' : \tld Y \rightarrow Y$ be log resolutions of $(X, \ba^c)$ and $(Y, (\ba \O_{\tld Y})^c)$ respectively.  Let $G$ be the divisor on $\tld X$ such that $\ba \O_{\tld X} = \O_{\tld X}(-G)$ and and $F$ be the divisor on $\tld Y$ such that $\ba \O_{\tld Y} = \O_{\tld Y}(-F)$.  We can choose these resolutions so that there is  map $\gamma : \tld Y \rightarrow \tld X$ such that the following diagram commutes. Note $\gamma^* G = F$.
\[
\xymatrix{
 \tld Y \ar[r]^{\gamma} \ar[d]_{\pi'} & \tld X \ar[d]^{\pi'} \\
 Y \ar[r]_{\rho} & X \\
}
\]
We will show that there is a natural map
\[
\O_{\tld X}(\lfloor cG \rfloor) \rightarrow \gamma_* \O_{\tld Y}(\lfloor cF \rfloor) = \gamma_* \O_{\tld Y}(\lfloor c(\gamma^*G) \rfloor).
\]
By composition with the map
\[
\O_{\tld X}(\lfloor cG \rfloor) \rightarrow \gamma_* \gamma^* \O_{\tld X}(\lfloor cG \rfloor) = \gamma_* \O_{\tld Y}(\gamma^*\lfloor cG \rfloor),
\]
we see it is sufficient to show that there is a natural inclusion $\O_{\tld Y}(\gamma^* \lfloor cG \rfloor) \subseteq \O_{\tld Y}(\lfloor c (\gamma^* G) \rfloor)$.  But this is true because even though round-down does not commute with pullbacks, there is always an inequality, $\gamma^* \lfloor cG \rfloor \leq \lfloor c (\gamma^* G)\rfloor$.

Now consider the following diagram,
\[
\xymatrix{
\myR \rho_* \myR \pi'_* \O_{\tld Y}(\lfloor c F \rfloor) & \myR \pi_* \O_{\tld X}(\lfloor c G \rfloor) \ar[l] \\
\myR \rho_* \O_Y \ar[u]_{\myR \rho_* (p') } & \O_X \ar[u]^{p} \ar[l] \\
}
\]
Since $(Y, (\ba \O_Y)^c)$ is rational, $p' : \O_Y \rightarrow \myR \pi_*' \O_{\tld Y}(\lfloor cF\rfloor)$ is a quasi-isomorphism.  Thus, consider the composition
\[
\O_X \rightarrow \myR \pi_* \O_{\tld X}(\lfloor c G \rfloor) \rightarrow \myR \rho_* \myR \pi'_* \O_{\tld Y}(\lfloor c F \rfloor) \qis \myR \rho_* \O_Y \rightarrow \O_X
\]
where the final map in the composition exists by hypothesis.  This composition must be a quasi-isomorphism by construction, creating a left inverse of $p$.  This completes the proof by Theorem \ref{GeneralizedKovacsSplittingTheorem}.
\end{proof}

\begin{cor}
\label{rationalSummands}
Suppose $R$ and $S$ are domains, $\ba$ is an ideal of $R$ and $R$ is a summand of $S$ (for example, suppose that $R$ is normal and that $R \rightarrow S$ is a finite map).  If $(S, (\ba S)^c)$ is rational, so is $(R, \ba^c)$.
\end{cor}

Compare the Corollary \ref{rationalSummands} with Proposition \ref{SummandTheoremForPositiveCharacteristic}.

\begin{rem}
Note that the converse to this statement is not true.  Of course, even when $\ba = R$, the converse can fail for a canonical cover by \cite{SinghCyclicCoversOfRational}.  We give an example another type of failure using the Example \ref{RationalNotKLTExample}.  Let $X = \Spec \bC[x^3, x^2 y, x y^2, y^3]$ and let $Y = \bC[x,y]$ be its canonical cover.  Let $\ba = (x^3, x^2 y, x y^2, y^3)$ and note that $\ba \O_Y = (x,y)^3$.  Then $(X, \ba^{0.9})$ has rational singularities, but $(Y, ((x,y)^3)^{0.9})$ clearly does not.
\end{rem}

We now explore how rational pairs deform; see \cite[Theorem 2]{ElkikDeformationsOfRational} or \cite[11.15]{KollarSingularitiesOfPairs}

\begin{thm}
\label{rationalDeforms}
Suppose that $(X, \ba^c)$ is a pair, that $H$ is a Cartier divisor on $X$, and that $H$ has no common components with $V(\ba)$.  If the pair $(H, (\ba|_H)^c)$ has rational singularities, then so does $(X, \ba^c)$ near $H$.
\end{thm}
\begin{proof} \cite{ElkikDeformationsOfRational}
Let $x$ be a point of $X$ also contained in $H$.  Note it is enough to prove the problem at the stalk associated to $x$ and so we assume that $X = \Spec R$ with $(R, m)$ local, $H = \Spec R / f$ for some regular element $f \in R$ and that $\ba$ is an ideal of $R$ which has no common minimal primes with $(f) = I_H$.  Note that since $(H, (\ba|_H)^c)$ is rational, $H$ and thus $X$ is Cohen-Macaulay.  Let $\pi : \tld X \rightarrow X$ be a resolution of $(X, \ba^c)$ that is also simultaneously a resolution of $H$ and let $G$ denote the divisor such that $\ba \O_{\tld X} = \O_{\tld X}(-G)$.  Let $\overline H$ be the total transform of $H$ (that is, $\overline H$ is the scheme defined by $f \O_{\tld X}$) and let $\tld H$ denote the strict transform of $H$.  Note, there is a natural inclusion of schemes $\tld H \rightarrow \overline H$.  Consider the following diagram.
\[
\scriptsize
\xymatrix{
& & &  \pi_* (\omega_{\tld H} \tensor \O_{\tld X}(-\lfloor c G\rfloor)) \ar[d]  \\
0 \ar[r] & \pi_* (\omega_{\tld X} \tensor \O_{\tld X}(-\lfloor c G\rfloor)) \ar[r]^-{\times f} \ar@{^{(}->}[d]^{\psi} & \pi_* (\omega_{\tld X}\tensor \O_{\tld X}(-\lfloor c G\rfloor)) \ar[r] \ar@{^{(}->}[d]^{\psi} & \pi_* (\omega_{\overline H} \tensor \O_{\tld X}(-\lfloor c G\rfloor)) \ar[r] \ar@{->>}[d]^{\phi} & 0 \\
0 \ar[r] & \omega_X \ar[r]^-{\times f} & \omega_X \ar[r] & \omega_H \ar[r] & 0
}
\]
The bottom row is exact because $H$ is Cohen-Macaulay.  The top row is exact by \ref{LocalVanishingForMultiplierSubmodules}.  The map labeled $\phi$ is surjective since the vertical composition from $\pi_* (\omega_{\tld H} \tensor \O_{\tld X}(-\lfloor c G\rfloor))$ is an isomorphism.  It is then enough to show that $\psi$ is surjective by \ref{LocalVanishingForMultiplierSubmodules}.

Let $C$ be the cokernel of $\psi$.  The fact that $\phi$ is surjective means that $\xymatrix{C \ar[r]^-{\times f} & C}$ is surjective by the snake lemma.  But this contradicts Nakayama's lemma, completing the proof.
\end{proof}

We conclude this section with several results related to adjunction.  Compare these results with \cite[9.5.11, 9.5.17]{LazarsfeldPositivity2}, \cite[5.6]{KollarMori}, \cite{ShokurovThreeDimensionalLogFlips}.  The first result could be thought of as an analogue of adjunction and inversion of adjunction for log terminal singularities, and in some sense it is the easy case, since we work only with Cartier divisors; compare with \cite[Chapters 16 and 17]{KollarFlipsAndAbundance}.  We also obtain a positive characteristic analogue later in Theorem \ref{inversion}.

\begin{thm}
\label{purelyRationalInversionOfAdjunction}
Suppose that $X$ is a normal scheme and that $H$ is a Cartier divisor on $X$.  Further suppose that $\ba$ is an ideal sheaf whose support does not contain any component of $H$ and $c$ is a nonnegative real number.  Then $(H, (\ba|_H)^c)$ has rational singularities if and only if $(X, H;\ba^c)$ has purely rational singularities near $H$.
\end{thm}

\begin{proof}
By remark \ref{adjointMultiplierShortExactSequence}, we have a short exact sequence which maps to another short exact sequence
\[
\xymatrix{
0 \ar[r] & \mJ(\omega_X, \ba^c) \ar[r] \ar[d]^{\alpha} &  \adj(\omega_X, H; \ba^c) \tensor \O_X(H) \ar[r] \ar[d]^{\beta} &  \mJ(\omega_H, (\ba|_H)^c) \ar[r] \ar[d]^{\gamma} &  0 \\
0 \ar[r] & \omega_X \ar[r] &  \omega_X(H) \ar[r] &  \omega_H \ar[r] &  0.
}
\]
The bottom row is exact on the right since $X$ is Cohen-Macaulay near $H$ under any assumption.

Suppose first that $(H, (\ba|_H)^c)$ has rational singularities, then so does $(X, \ba^c)$ near $H$.  By localizing, we assume that $(X, \ba^c)$ is rational.  These observations imply that the maps $\alpha$ and $\gamma$ are isomorphisms which proves that $\beta$ is an isomorphism as well.  Untwisting by $\O_X(H)$ implies that $(X, H;\ba^c)$ has purely rational singularities.

Conversely, if $(X, H;\ba^c)$ has purely rational singularities then $(X, \ba^c)$ has rational singularities by remark \ref{PurelyRationalImpliesRationalRemark}.  Thus $\alpha$ and $\beta$ are isomorphisms which implies that $\gamma$ is an isomorphism as well, which completes the proof.
\end{proof}

\begin{rem}
One could of course dualize the proof of \ref{purelyRationalInversionOfAdjunction} and perform the same argument in the derived category.  In the case that $H$ wasn't Cartier, one could prove the same result using the suggested definition from Remark \ref{RemarkOnNonCartierCase} if one assumed that $(X, \ba^c)$ already had rational singularities for the rational implies purely rational (``only if'') implication.
\end{rem}

We also have the following result which can be thought of as an analogue to the ``adjunction direction'' for log canonical singularities.

\begin{thm}
\label{rationalDuBoisInversionOfAdjunction}
Suppose that $X$ is a reduced scheme and that $H$ is a Cartier divisor.  If $(X, (1-\epsilon)H)$ is rational for all sufficiently small $\epsilon > 0$, then $H$ has Du Bois singularities.
\end{thm}
\begin{proof}
Let $\pi : \tld X \rightarrow X$ be a log resolution of $(X, H)$.  Let $\overline H$ be the total transform of $H$ and let $E$ be the reduced pre-image of $H$ under $\pi$, in particular $\overline{H}_{\red} = E$.  Consider the following diagram.  Note that since $\O_X \qis \myR \pi_* \O_{\tld X}(\lfloor (1-\epsilon)\overline{H}\rfloor)$ for all $\epsilon$ sufficiently close to zero, we have
\[
\O_X(-H) \qis \myR \pi_* \O_{\tld X}(\lfloor(1-\epsilon)\overline{H} - \overline H \rfloor) \qis \myR \pi_* \O_{\tld X}(\lfloor -\epsilon\overline{H} \rfloor) \qis \myR \pi_* \O_{\tld X}(-E)
\]
for $\epsilon$ sufficiently small.  Therefore
\[
\xymatrix{
 \O_{X}(-H) \ar[d] \ar[r] & \O_{X} \ar[r] \ar[d] & \O_H \ar[d] \ar[r]^{+1} &  \\
\myR \pi_* \O_{\tld X}(-E) \ar[r] & \myR \pi_* \O_{\tld X} \ar[r] & \myR \pi_* \O_{E} \ar[r]^-{+1} &\\
}
\]
and the first two vertical arrows are quasi-isomorphisms.  But then the third arrow is also a quasi-isomorphism which proves that $H$ has Du Bois singularities by \cite[2.4]{KovacsDuBoisLC1}; also see \cite[12.8]{KollarShafarevich}.
\end{proof}

There is a partial converse to Theorem \ref{rationalDuBoisInversionOfAdjunction}, which can be thought of as an analogue to inversion of adjunction for log canonicity; compare with \cite{KawakitaInversion}.

\begin{thm}
\label{rationalDuBoisInversion2}
Suppose that $X$ is a reduced scheme and $H$ is a Cartier divisor on $X$.  Further suppose that $X \backslash H$ is smooth.  Then $H$ has Du Bois singularities if and only if $(X, (1-\epsilon)H)$ is rational near $H$ for all sufficiently small $\epsilon > 0$.
\end{thm}
\begin{proof}
We setup the proof in the same way as \ref{rationalDuBoisInversionOfAdjunction}, but now we note that $\O_H \qis \myR \pi_* \O_{E}$ if and only if $H$ has Du Bois singularities since $X-H$ is smooth.  Note that if $H$ is Du Bois, $X$ automatically has rational singularities (and thus is also Cohen-Macaulay) by \cite[5.1]{SchwedeEasyCharacterization}.
\end{proof}

\begin{rem}
When working with any ambient $X$ with rational singularities, see \cite{SchwedeEasyCharacterization}, we expect that $\O_H \qis \myR \pi_* \O_{E}$ if and only if $H$ has Du Bois singularities.  Therefore, the hypothesis that $X-H$ is smooth could possibly be replaced with the condition that $X-H$ is rational without otherwise altering the proof.
\end{rem}

\begin{rem}
We also have positive characteristic analogue of the previous two theorems using $F$-injective instead of Du Bois singularities, see Proposition \ref{FInjectiveInversion}.
\end{rem}

It is a conjecture of Koll\'ar that log canonical singularities are Du Bois and the previous proof shows that this conjecture is closely related to inversion of adjunction on log canonicity.  Recent work by Kov\'acs, the first author and Smith has given a proof that (semi)log canonical singularities are Du Bois in the case of Cohen-Macaulay schemes; see \cite{KovacsSchwedeSmithLCImpliesDuBois}.  That result and the previous argument also give a very short homological proof of inversion of adjunction on log canonicity (at least in the case that $X$ has Gorenstein rational singularities and is smooth outside $H$).

The previous result suggests that it might be natural to consider Du Bois singularities for pairs, and perhaps even suggests a definition.  However, there are certain technicalities associated to such a definition when the ambient space is not ``nice''.
In positive characteristic, we do propose an analogous definition, at least in the Cohen-Macaulay case; see Definition \ref{DefinitionFInjectiveThreshold}.

\section{Positive characteristic preliminaries}
\label{SectionPositiveCharacteristicPreliminaries}
In this section, we recall the definitions of generalizations of tight closure and $F$-singularities of pairs.
The reader is referred to \cite{HaraYoshidaGeneralizationOfTightClosure}, \cite{TakagiInversion}, \cite{TakagiPLTAdjoint}, \cite{TakagiWatanabeFPureThresh} and \cite{TakagiYoshidaGeneralizedTestIdealsAndSymbolicPowers} for details.

Throughout the following sections, all rings are excellent reduced Noetherian commutative rings with identity.
Let $R$ be a reduced ring of characteristic $p>0$.
We denote by $R^{\circ}$ the set of elements of $R$ that are not in any minimal prime ideal of $R$.
Let $F:R \to R$ be the Frobenius map which sends $x$ to $x^p$.
$R$ viewed as an $R$-module via the $e$-times iterated Frobenius map $F^e \colon R \to R$ is denoted by ${}^e\! R$.
Since $R$ is reduced, we can identify $F^e:R \to {}^e\! R$ with the natural inclusion map $R \hookrightarrow R^{1/p^e}$.
Also, for any ideal $I$ of $R$ and for any power $q$ of $p$, we denote by $I^{[q]}$ the ideal of $R$ generated by the $q$-th powers of elements of $I$. We say that $R$ is {\it $F$-finite} if ${}^1\! R$ (or $R^{1/p}$) is a finitely generated $R$-module.
For example, any algebra essentially of finite type over a perfect field is $F$-finite.

Let $M$ be an $R$-module.
For each integer $e \ge 1$, we denote $\F^{e}(M)= \F_{R}^{e}(M) := {}^e\! R \otimes_R M$ and
regard it as an $R$-module by the action of $R$ from the left.
Then we have the induced $e$-times  Frobenius map $F^{e} \colon M \to
\F^{e}(M)$. The image of $z \in M$ via this map is denoted by $z^q:= F^{e}(z)
\in \F^{e}(M)$.
For an $R$-submodule $N$ of $M$, we denote by $N^{[q]}_{M}$ the
image of the induced map $\F^{e}(N) \to \F^{e}(M)$.
If $M=R$ and $N$ is an ideal $I$ of $R$, then $I^{[q]}_{R}=I^{[q]}$.

\begin{defn}[\textup{\cite[Definition 6.1]{HaraYoshidaGeneralizationOfTightClosure}, \cite[Definition 3.1]{TakagiPLTAdjoint}}]\label{tight closure}
Let $R$ be a reduced ring of characteristic $p>0$, $\a \subseteq R$ be an ideal such that $\a \cap R^{\circ} \ne \emptyset$ and $t \ge 0$ be a real number.
Let $N \subseteq M$ be (not necessarily finitely generated) $R$-modules.
\renewcommand{\labelenumi}{(\roman{enumi})}
\begin{enumerate}
\item
The \textit{$\a^t$-tight closure} $N^{*\a^t}_M$
of $N$ in $M$ is defined to be the submodule of $M$ consisting of all elements $z \in M$ for which there exists $c \in R^{\circ}$ such that
$$c\a^{\lceil tq \rceil}z^q \subseteq N^{[q]}_M$$
for all large $q = p^e$. The $\a^t$-tight closure of an ideal $I$ of $R$ is simply defined by $I^{*\a^t}:=I_R^{*\a^t}$.
\item
Let $x \in R^{\circ}$ such that $\a$ is not contained in any minimal prime of $xR$.
Then the \textit{divisorial $(x;\a^t)$-tight closure} $N^{\mathrm{div}*(x;\a^t)}_M$
of $N$ in $M$ is defined to be the submodule of $M$ satisfying the following condition:
an element $z \in M$ belongs to $N^{\mathrm{div}*(x;\a^t)}_M$ if there exists $c \in R^{\circ}$ which is not in any minimal prime of $xR$ such that
$$cx^{q-1}\a^{\lceil tq \rceil}z^q \subseteq N^{[q]}_M$$
for all large $q = p^e$. The divisorial $(x;\a^t)$-tight closure of an ideal $I$ of $R$ is simply defined by $I^{\mathrm{div}*(x;\a^t)}:=I_R^{\mathrm{div}*(x;\a^t)}$.
\end{enumerate}
\end{defn}

\begin{rem}
When $\a=R$, $\a^t$-tight closure is nothing but classical tight closure, that is, the classical tight closure $I^*$ of an ideal $I \subseteq R$ is equal to $I^{*R^t}$ for any $t \ge 0$.
We refer the reader to \cite{HunekeTightClosureBook} for the classical tight closure theory.
\end{rem}

\begin{defn}[\textup{\cite[Definition 6.3]{HaraYoshidaGeneralizationOfTightClosure}}]
Let $R, \a, t$ be as in Definition \ref{tight closure}.
An element $c \in R^{\circ}$ is called an {\it $\a^t$-test element} if for every ideal $I \subseteq R$, we have $cz^q\a^{\lceil tq \rceil} \subseteq I^{[q]}$ for all $q = p^e$ whenever $z \in I^{*\a^t}$.
\end{defn}

A local ring $R$ of characteristic $p>0$ is said to be \textit{$F$-rational} if $I^*=I$ for all ideals $I \subseteq R$ generated by a system of parameters for $R$ (see \cite{FedderWatanabe} for details).
\begin{lem}[\textup{\cite{TakagiYoshidaGeneralizedTestIdealsAndSymbolicPowers}}]\label{completion}
Let $(R,\m)$ be an excellent reduced local ring of characteristic $p>0$.
Let $\a \subseteq R$ be an ideal such that $\a \cap R^{\circ} \ne \emptyset$ and $t \ge 0$ be a real number.
 \begin{enumerate}
 \item Let $\widehat{R}$ denotes the $\m$-adic completion of $R$. Then $I^{*\a^t}\widehat{R}=(I\widehat{R})^{*(\a\widehat{R})^t}$ for all $\m$-primary ideals $I$ of $R$.

\item If $R$ is equidimensional and $S$ is a multiplicatively closed set in $R$, then $I^{*\a^t}R_S=(IR_S)^{*(\a R_S)^t}$ for all ideals $I$ generated by a subsystem of parameters for $R$.

\item Let $c \in R^{\circ}$ such that $R_c$ is Gorenstein $F$-rational. Then some power $c^n$ of $c$ is an $\a^t$-test element for all ideals $\a \subseteq R$ such that $\a \cap R^{\circ} \ne \emptyset$ and for all real numbers $t \ge 0$.
\end{enumerate}
\end{lem}

\begin{defn}[\textup{\cite[Definition 3.1]{TakagiInversion}}]\label{F-def}
Let $\a$ be an ideal of an $F$-finite reduced ring $R$ of characteristic $p>0$ such that $\a \cap R^{\circ} \ne \emptyset$ and let $t \ge 0$ be a real number.
\renewcommand{\labelenumi}{(\roman{enumi})}
\begin{enumerate}
\item The pair $(R, \a^t)$ is said to be {\it $F$-pure} if for all large $q=p^e$, there exists an element $d \in \a^{\lfloor t(q-1) \rfloor}$ such that the natural inclusion $d^{1/q}R \hookrightarrow R^{1/q}$ splits as an $R$-module homomorphism.
\item The pair $(R, \a^t)$ is said to be {\it strongly $F$-regular} if for every $c \in R^{\circ}$, there exist $q=p^e$ and $d \in \a^{\lceil tq \rceil}$ such that the natural inclusion $(cd)^{1/q}R \hookrightarrow R^{1/q}$ splits as an $R$-module homomorphism.
\item Let $x \in R^{\circ}$ such that $\a$ is not contained in any minimal prime of $xR$.  The triple $(R,x;\a^t)$ is said to be \textit{divisorially $F$-regular} if for every $c \in R^{\circ}$ which is not in any minimal prime of $xR$, there exist $q=p^e$ and $d \in \a^{\lceil tq \rceil}$ such that the natural inclusion $(cdx^{q-1})^{1/q}R \hookrightarrow R^{1/q}$ splits as an $R$-module homomorphism.
\end{enumerate}
\end{defn}

\begin{defn}[\textup{\cite[Definition 2.1]{TakagiWatanabeFPureThresh}}]
Let $R, \a$ be as in Definition \ref{F-def}. Assume in addition that $R$ is a strongly $F$-regular ring, that is, the pair $(R,R^1)$ is strongly $F$-regular.
Then the $F$-pure threshold $\mathrm{fpt}(\a)$ of $\a$ is defined to be
\begin{align*}
\mathrm{fpt}(\a)&=\{t \in \R_{\geq 0} \mid (R,\a^t)\textup{ is strongly $F$-regular}\}\\
&=\{t \in \R_{\geq 0} \mid (R,\a^t)\textup{ is $F$-pure}\}.
\end{align*}
\end{defn}

\begin{rem}\label{F-pair remark}
(1) When $\a=R$, the strong $F$-regularity (resp. $F$-purity) of $(R,\a^t)$ is equivalent to that of $R$.
 We refer the reader to \cite{HochsterHunekeTC1}, \cite{HochsterHunekeTightClosureAndStrongFRegularity} and \cite{HochsterRobertsFrobeniusLocalCohomology} for $F$-pure rings and strongly $F$-regular rings.

(2) If $(R,\a^t)$ is strongly $F$-regular, then it is $F$-pure. If $(R,x;\a^t)$ is divisorially $F$-regular, then $(R,x\a^t)$ is $F$-pure and $(R,x^{1-\epsilon}\a^t)$ is strongly $F$-regular for any $1 \ge \epsilon >0$.
The reader is referred to \cite{HaraWatanabeFRegFPure}.

(3) If $(R,\a^t)$ is strongly $F$-regular (resp. $(R,x;\a^t)$ is divisorially $F$-regular), then $I^{*\a^t}=I$ (resp. $I^{\mathrm{div}*(x,\a^t)}=I$) for all ideals $I \subseteq R$.
If $R$ is $F$-finite $\Q$-Gorenstein, then the converse also holds true. The reader is referred to \cite[Corollary 3.5]{TakagiInversion} (resp. \cite[Remark 3.2]{TakagiPLTAdjoint}).
\end{rem}

\section{Basic definitions and fundamental properties in positive characteristic}
\label{SectionDefinitionsAndPropertiesInPositiveCharacteristic}
In \cite{FedderWatanabe}, Fedder and Watanabe defined the notion of $F$-rational rings.
In this section, we introduce the notion of $F$-rationality for a pair $(R,\a^t)$ of a ring $R$ of characteristic $p>0$  and an ideal $\a \subseteq R$ with real exponent $t \ge 0$.

\begin{defn}[cf.\textup{\cite{FedderWatanabe}}]\label{F-pair}
Let $\a$ be an ideal of  a reduced ring $R$ of characteristic $p>0$ such that $\a \cap R^{\circ} \ne \emptyset$ and let $t \ge 0$ be a real number.

When $R$ is local, $(R,\a^t)$ is said to be \textit{$F$-rational} if $I^{*\a^t}=I$ 
for every ideal $I$ generated by  a system of parameters for $R$.
When $R$ is not local, we say that $(R,\a^t)$ is $F$-rational if the localization $(R_{\m},\a_{\m}^t)$ is $F$-rational for every maximal ideal $\m$ of $R$.
\end{defn}

\begin{prop}\label{basic}
Let $\a \subseteq \b$ be ideals of a reduced ring $R$ of characteristic $p>0$ such that $\a \cap R^{\circ} \ne \emptyset$ and let $t \ge 0$ be a real number.
\begin{enumerate}
\item If $(R,\a^t)$ is $F$-rational, then so is $(R,\a^{s})$ for all $0 \le s \le t$.
\item If $(R,\a^t)$ is $F$-rational, then so is $(R,\b^{t})$.
When $\a$ is a reduction of $\b$, $(R,\a^t)$ is $F$-rational if and only if $(R,\b^t)$ is $F$-rational.
\item If $(R,\a^t)$ is $F$-rational, then $R$ is $F$-rational and is, in particular, normal.
Moreover, if $R$ is locally excellent, then $R$ is Cohen-Macaulay.
\end{enumerate}
\end{prop}
\begin{proof}
It is immediate from \cite[Proposition 1.3]{HaraYoshidaGeneralizationOfTightClosure} and \cite[Theorem 4.2]{HunekeTightClosureBook}.
\end{proof}

\begin{lem}\label{local characterization}
Let $(R,\m)$ be a $d$-dimensional excellent reduced local ring of characteristic $p>0$.
Let $\a$ be an ideal of $R$ such that $\a \cap R^{\circ} \ne \emptyset$ and let $t \ge 0$ be a real number.
Then the following three conditions are equivalent to each other.
\begin{itemize}
\item[(1)] $(R,\a^t)$ is $F$-rational.
\item[(2)] $R$ is equidimensional and $I^{*\a^t}=I$ for some ideal $I$ generated by a system of parameters for $R$.
\item[(3)]  $R$ is Cohen-Macaulay and $0^{*\a^t}_{H^d_{\m}(R)}=0$ in $H^d_{\m}(R)$.
\end{itemize}
Furthermore, if $\ba = (f)$ is principal, then the previous three conditions are equivalent to the following:
\begin{itemize}
\item[(4)] $R$ is Cohen-Macaulay and for each $c \in R^{\circ}$, there exists $q=p^e$ such that $cf^{\lceil tq \rceil} F^e:H^d_{\m}(R) \to H^d_{\m}(R)$ is injective.
\end{itemize}
\end{lem}

\begin{proof}
The implication (1) $\Rightarrow$ (2) is obvious.
So, we will prove the implication (2) $\Rightarrow$ (3).
First note that $R$ is Cohen-Macaulay, because $I^{*\a^t}=I^*=I$ and $R$ is equidimensional (see \cite[Theorem 4.2]{HunekeTightClosureBook}).
We choose a system of parameters $x_1, \dots, x_d$ in $R$ such that $I=(x_1, \dots, x_d)$ and let $x$ denote the product of $x_1, \dots, x_d$.

\begin{cl}
$(x_1^m, \dots, x_d^m)^{*\a^t}=(x_1^m, \dots, x_d^m)$ for each integer $m \ge 1$.
\end{cl}
\begin{proof}[Proof of Claim]
Let $y \in (x_1^m, \dots, x_d^m)^{*\a^t}$.
Without loss of generality we may assume that $y(x_1, \dots, x_d) \subseteq (x_1^m, \dots, x_d^m)$.
Since $R$ is Cohen-Macaulay, one has $y \in (x_1^m, \dots, x_d^m, x^{m-1})$.
We write down $y$ as $y=\sum_{i=1}^d a_ix_i^m + bx^{m-1}$ where $a_i \in R$ for all $i=1, \dots, d$ and $b \in R$.
By definition, there exists $c \in R^{\circ}$ such that $c\a^{\lceil tq \rceil}y^q \in (x_1^{mq}, \dots, x_d^{mq})$ for all large $q=p^e$.
Then $c\a^{\lceil tq \rceil}b^q \in (x_1^q, \dots, x_d^q)$.
Hence $b \in (x_1, \dots, x_d)^{*\a^t}=(x_1, \dots, x_d)$, which implies that $y \in (x_1^m, \dots, x_d^m)$.
\end{proof}

Fix an arbitrary element $\eta=[\frac{z}{x^m}] \in 0^{*\a^t}_{H^d_{\m}(R)}$.
By the definition of $\a^t$-tight closure, there exists $c \in R^{\circ}$ such that $0=c\a^{\lceil tq \rceil}\eta^q=c\a^{\lceil tq \rceil}[\frac{z^q}{x^{mq}}]$ for all large $q=p^e$.
This implies that for large $n$, $c\a^{\lceil tq \rceil}z^qx^n \in (x_1^{n+mq}, \dots, x_d^{n+mq})$.
As $R$ is Cohen-Macaulay, we then obtain that $c\a^{\lceil tq \rceil}z^q \in (x_1^{mq}, \dots, x_d^{mq})$ for all large $q=p^e$ which gives that $z \in (x_1^m, \dots, x_d^m)^{*\a^t}=(x_1^m, \dots, x_d^m)$, where the last equality follows from the above claim.
Then $\eta=0$, that is,  $0^{*\a^t}_{H^d_{\m}(R)}=0$.

Next we will show the implication (3) $\Rightarrow$ (1).
Take any system of parameters $x_1, \dots, x_d$ in $R$ and let $x$ represent the product of $x_1, \dots, x_d$.
Fix any element $z \in (x_1, \dots, x_d)^{*\a^t}$ and consider the element $\xi=[\frac{z}{x}] \in H^d_{\m}(R)$.
By definition, we can choose an element $d \in R^{\circ}$ such that $d\a^{\lceil tq \rceil}z^q \in (x_1^q, \dots, x_d^q)$ for all large $q=p^e$.
This implies that $d\a^{\lceil tq \rceil}\xi^q=0$ for all large $q=p^e$, that is, $\xi \in 0^{*\a^t}_{H^d_{\m}(R)}$.
By assumption, one has $\xi=0$ which gives $z \in (x_1, \dots, x_d)$.

Finally we note that the implication (4) $\Rightarrow$ (3) is obvious, so it suffices to show that (3) $\Rightarrow$ (4).  Therefore, we assume that $0^{*f^t}_{H^d_{\m}(R)} = 0$.  Notice that for each element $c \in
R^{\circ}$, and non-zero $z \in H^d_{\m}(R)$, there exists $e > 0$ such that $c f^{\lceil t q \rceil} z^q \neq 0$ (where $q = p^e$).  It then follows for every $q'$, since the ring itself is $F$-injective, that
\[c^{q'} f^{q' \lceil t q \rceil} z^{qq'} \neq 0 \text{ which implies that }
c f^{\lceil t qq' \rceil} z^{qq'} \neq 0.\]
This implies that for all sufficiently large $e$, that $c f^{\lceil t q \rceil} z^q \neq 0$.

Fix a $c \in R^{\circ}$ and consider the modules $N_e \subset H^d_{\m}(R)$ defined as
\[
N_e = \text{ker} (c f^{\lceil t p^e \rceil} F^e : H^d_{\m}(R) \rightarrow H^d_{\m}(R))
\]
The previous work guarantees that these modules form a decreasing sequence in $H^d_{\m}(R)$, an Artinian module.  On the other hand, no non-zero element is in the intersection of all of the $N_e$.  Therefore, they must stabilize at zero at some finite step.  This implies condition (4).
\end{proof}

\begin{rem}
\label{testModuleDefinition}
Let the notation be as in Lemma \ref{local characterization} and assume in addition that $R$ is a homomorphic image of a Gorenstein local ring.
Then one could define the \textit{generalized parameter test submodule} $\tau(\omega_R,\a^t)$ associated to $(R,\a^t)$ to be
$$\tau(\omega_R,\a^t)=\mathrm{Ann}_{\omega_R}(0^{*\a^t}_{H^d_{\m}(R)}) \subseteq \omega_R.$$
This is a characteristic $p$ analogue of the multiplier submodule (see Definition \ref{MultiplierSubmoduleDefinition}).
It follows that $(R,\a^t)$ is $F$-rational if and only if $R$ is Cohen-Macaulay and $\tau(\omega_R,\a^t)=\omega_R$.
Employing the same strategy as that of \cite{HyryVillamayorBriansconSkodaForIsolated}, we can use the generalized parameter test submodule $\tau(\omega_R,\a^t)$ to recover the Brian\c con-Skoda theorem for $F$-rational rings (\cite[Theorem 3.6]{AberbachHunekeFRationalAndIntegralClosures}):  if $(R,\m)$ is an excellent $F$-rational local ring of dimension $d$ which is a homomorphic image of a Gorenstein local ring, then $\overline{I^{n+d-1}} \subseteq I^n$ for all ideals $I \subseteq R$ and integers $n \ge 0$.
\end{rem}

\begin{prop}\label{basic2}
Let $\a$ be an ideal of a locally excellent reduced ring $R$ of characteristic $p>0$ such that $\a \cap R^{\circ} \ne \emptyset$ and let $t \ge 0$ be a real number.
\begin{enumerate}
\item If $(R,\a^t)$ is strongly $F$-regular, then it is $F$-rational.
If $R$ is $F$-finite Gorenstein, then the converse also holds true.
\item Let $S$ be a multiplicatively closed set in $R$.
If $(R,\a^t)$ is $F$-rational, then the localization $(R_S, \a_S^t)$ is also $F$-rational.
\item Assume in addition that $R$ is local.
Then $(R,\a^t)$ is $F$-rational if and only if $(\widehat{R}, (\a\widehat{R})^t)$ is $F$-rational.
\end{enumerate}
\end{prop}
\begin{proof}
(1)
By Remark \ref{F-pair remark} (3),
strongly $F$-regular pairs are $F$-rational. So, we consider the converse implication.
Since strong $F$-regularity commutes with localization, we may assume that $(R,\m)$ is an $F$-finite reduced local ring.
By \cite[Corollary 3.5]{TakagiInversion},
$(R,\a^t)$ is strongly $F$-regular if and only if $0^{*\a^t}_E=0$, where $E=E_R(R/\m)$ is the injective hull of the residue field $R/\m$.
Thus,
if $R$ is Gorenstein,
then by Lemma \ref{local characterization},
the F-rationality of $(R,\a^t)$ is equivalent to the strong F-regularity of $(R,\a^t)$
since in this case $H^{\dim R}_{\m}(R) \cong E$.

(2) We may assume that $R$ is a Cohen-Macaulay local ring, and it suffices to show that $(R_P, \a_P^t)$ is $F$-rational for every prime ideal $P$ of $R$.
Let $x_1, \dots, x_i$ be any elements of $P$ whose images in $R_P$ form a system of parameters for $R_P$.
We can choose elements $x_{i+1}, \dots, x_d$ of $R$ such that $x_1, \dots, x_d$ form a system of parameters for $R$.
Set $I=(x_1, \dots, x_i)$ and $I_n=(x_1, \dots, x_i, x_{i+1}^n, \dots, x_d^n)$ for each integer $n \ge 1$.
By assumption, one has $I_n^{*\a^t}=I_n$ for all $n \ge 1$.
This implies that
$$I=\bigcap_n I_n=\bigcap_n I_n^{*\a^t}=I^{*\a^t}.$$
Since $I$ is generated by a subsystem of parameters for $R$, by Lemma \ref{completion} (2), one has $(IR_P)^{*\a_P^t}=I^{*\a^t}R_P=IR_P$.
That is, $(R_P,\a_P^t)$ is $F$-rational.

(3) Let $I$ be an ideal of $R$ generated by a system of parameters for $R$.
By Lemma \ref{completion} (1), $I^{*\a^t}=I$ if and only if $(I\widehat{R})^{*(\a\widehat{R})^t}=I \widehat{R}$.
Thus, the assertion is obvious.
\end{proof}

\begin{defn}\label{test def}
Let $\a$ be an ideal of a reduced local ring $R$ of characteristic $p>0$ such that $\a \cap R^{\circ} \ne \emptyset$ and let $t \ge 0$ be a real number.
Then an element $c \in R^{\circ}$ is called a {\it parameter $\a^t$-test element} if for every ideal $I$ generated by a system of parameters for $R$, we have $cz^q\a^{\lceil tq \rceil} \subseteq I^{[q]}$ for all $q=p^e$ whenever $z \in I^{*\a^t}$.
\end{defn}

\begin{rem}\label{test rem}
Let $R$ be a Cohen-Macaulay reduced local ring of characteristic $p>0$ and $c \in R^{\circ}$ be a parameter $\a^t$-test element for a principal ideal $\a = (f)$.
Then, by the same argument as the proof of Lemma \ref{local characterization}, we can easily check that $(R,\a^t)$ is $F$-rational if and only if there exist $q=p^e$ and $c' \in \a^{\lceil tq \rceil}$ such that $cc' F^e:H^d_{\m}(R) \to H^d_{\m}(R)$ is injective,
where $F^e:H^d_{\m}(R) \to H^d_{\m}(R)$ denotes the induced $e$-times iterated Frobenius map on $H^d_{\m}(R)$.
\end{rem}

\begin{lem}\label{test exist}
Let $(R,\m)$ be a $d$-dimensional excellent reduced equidimensional local ring of characteristic $p>0$.
Let $c \in R^{\circ}$ such that $R_c$ is $F$-rational. Then some power $c^n$ of $c$ is a parameter $\a^t$-test element for all ideals $\a \subseteq R$ such that $\a \cap R^{\circ} \ne \emptyset$ and for all real numbers $t \ge 0$.
\end{lem}
\begin{proof}
Making use of gamma construction, by an argument analogous to the proof of \cite[Theorem 3.9]{VelezOpennessOfTheFRationalLocus}, we can reduce to the case where $R$ is an $F$-finite reduced local ring which is a homomorphic image of a Gorenstein local ring.
Let $c' \in R^{\circ}$ be an $R$- and $\a^t$-test element (we can take such an element by Lemma \ref{completion} (3)) and let $F^e:H^d_{\m}(R) \to H^d_{\m}(R)$ denote the induced $e$-times iterated Frobenius map on $H^d_{\m}(R)$.

\begin{cl}[\textup{\cite[Theorem 1.13]{VelezOpennessOfTheFRationalLocus}}]
There exists $q_0=p^{e_0}$ and $n \in \N$ such that the $n^{\rm th}$ power $c^n$ of $c$ kills $\Ker(c'F^{e_0})$.
\end{cl}

Take any system of parameters $x_1, \dots, x_d$ in $R$ and let $x$ denote the product of $x_1 \dots x_d$.
Fix any $z \in (x_1, \dots, x_d)^{*\a^t}$ and consider the element $\xi=[\frac{z}{x}] \in H^d_{\m}(R)$.
Since $c'$ is an $\a^t$-test element, one has $c'\a^{\lceil tq_0q \rceil}z^{q_0q} \in (x_1^{q_0q}, \dots, x_d^{q_0q})$ for all $q=p^e$, which implies that $c'\a^{\lceil tq_0q \rceil}\xi^{q_0q}=0$ in $H^d_{\m}(R)$.
In particular, $\a^{\lceil t q \rceil}\xi^q$ is contained in $\Ker(c'F^{e_0})$ and, therefore, $c^n \a^{\lceil t q \rceil}\xi^q=0$ by the above claim.
Then there exists an integer $k \ge 0$ such that $c^n \a^{\lceil t q \rceil}z^qx^k \in (x_1^{q+k}, \dots, x_d^{q+k})$.
Applying the colon-capturing property of classical tight closure and \cite[Lemma 12.9]{HochsterAndHunekePhantom}, one has some power $c^m$ of $c$ such that $c^{m}c^n\a^{\lceil t q \rceil}z^q \in (x_1^{q}, \dots, x_d^{q})$ for all $q=p^e$.
Since $m$ is independent of the choice of $x_1, \dots, x_d, z, \a, t$, $c^{m+n}$ is an $\a^t$-test element for all ideals $\a \subseteq R$ such that $\a \cap R^{\circ} \ne \emptyset$ and for all real numbers $t \ge 0$.
\end{proof}

\begin{thm}[\textup{cf.~ \cite[Theorem 3.1]{SmithFRatImpliesRat}}]
Let $R$ be an excellent reduced local ring of characteristic $p>0$, let $\a$ be an ideal of $R$ such that $\a \cap R^{\circ} \ne \emptyset$ and let $t \ge 0$ be a real number.
If $(R,\a^t)$ is $F$-rational, then it is pseudo-rational.
\end{thm}
\begin{proof}
Since $R$ is excellent $F$-rational, it is Cohen-Macaulay, normal and analytically unramified.
Let $\pi$ and $\delta_{\pi}$ be as in Definition \ref{pseudo-rational}.
Then by \cite[Proposition 3.8]{HaraYoshidaGeneralizationOfTightClosure}, one has $\Ker(\delta_{\pi}) \subseteq 0^{*\a^t}_{H^d_{\m}(R)}$.
Since $(R,\a^t)$ is $F$-rational, by Lemma \ref{local characterization}, this implies that $\Ker(\delta_{\pi})=0$.
\end{proof}

Let $R$ be an algebra essentially of finite type over a field $k$ of characteristic zero.
Let $\a \subseteq  R$ be an ideal such that $\a \cap R^{\circ} \ne \emptyset$ and let $t \ge 0$ be a real number.
One can choose a finitely generated $\Z$-subalgebra $A$ of $k$ and a subalgebra $R_A$ of $R$ essentially of finite type over $A$ such that the natural map $R_A \otimes_A k \to R$ is an isomorphism and ${\a_A}R=\a$ where ${\a_A}:=\a \cap R_A \subseteq R_A$.
Given a closed point $s \in \Spec A$ with residue field $\kappa=\kappa(s)$, we denote the corresponding fibers over $s$ by $R_{\kappa}, {\a_{\kappa}}$.
Then we refer to a triple $(\kappa, R_{\kappa}, {\a_{\kappa}})$, for a general closed point $s \in \Spec A$ with residue field $\kappa=\kappa(s)$ of sufficiently large characteristic $p \gg 0$, as ``{\it reduction to characteristic $p \gg 0$}'' of $(k,R, \a)$.  The pair $(R_{\kappa}, {\a_{\kappa}}^{t})$ inherits the properties possessed by the original pair $(R, \a^{t})$ (the size of $p$ depends on $t$).
Furthermore, given a log resolution $f:\widetilde{X} \to X=\Spec R$ of $(X, \a)$, we can reduce this entire setup to characteristic $p \gg 0$.

\begin{defn}\label{reduction}
In the above situation,
$(R,\a^t)$ is said to be of \textit{strongly $F$-regular} (resp. \textit{$F$-pure, $F$-rational}) \textit{type} if the reduction to characteristic $p \gg 0$ of $(R,\a^t)$ is strongly $F$-regular (resp. $F$-pure, $F$-rational).
\end{defn}

\begin{thm}[\textup{cf.~\cite{HaraRatImpliesFRat}, \cite{MehtaSrinivasRatImpliesFRat}}]\label{correspondence}
Let $R$ be a finitely generated algebra over a field of characteristic zero.
Let $\a \subseteq R$ be an ideal such that $\a \cap R^{\circ} \ne \emptyset$ and $t \ge 0$ be a real number.
Then $(\Spec R,\a^t)$ has rational singularities if and only if $(R,\a^t)$ is of $F$-rational type.
\end{thm}
\begin{proof}

Since the assertion is local, we may assume that $(R,\m)$ is a $d$-dimensional normal Cohen-Macaulay local ring essentially of finite type over a field of characteristic zero.
Fix a log resolution $\pi:Y \to X:=\Spec R$ of $\a$ such that $\a\O_Y=\O_Y(-G)$, and let $E:=\pi^{-1}(\m)$ be the closed fiber of $\pi$ and let $\delta_{\pi}:H^d_{\m}(R) \to H^d_E(\O_Y(\lfloor tG \rfloor))$ be as in Definition \ref{pseudo-rational}.
Then by Remark \ref{pseudo rem}, $(\Spec R,\a^t)$ has rational singularities if and only if the map $\delta_{\pi}$ is injective.
After reduction to characteristic $p \gg 0$, we can assume that $R$ is a normal Cohen-Macaulay local ring essentially of finite type over a perfect field of characteristic $p$, together with a log resolution $\pi: Y \to X:=\Spec R$ of $(X,\a)$ such that  $\a \O_Y=\O_Y(-G)$.
Then it suffices to show that $(R,\a^t)$ is $F$-rational if and only if the map $\delta_{\pi}$ is injective, but it immediately follows from the combination of \cite[Theorem 6.9]{HaraYoshidaGeneralizationOfTightClosure} and Lemma \ref{local characterization}.
\end{proof}

\begin{rem}
In fact, using the same techniques, one can also give an equivalence between multiplier submodule and the parameter test submodule similar to \cite[Theorem 6.8]{HaraYoshidaGeneralizationOfTightClosure}.
\end{rem}

\begin{rem}
In \cite{SmithFRatImpliesRat}, Smith gave a characterization of $F$-rational rings in terms of the stability of submodules of $H^d_{\bm}(R)$ under the action of Frobenius.  Using the technique of Hara and Yoshida, see \cite[Proposition 1.15]{HaraYoshidaGeneralizationOfTightClosure}, one can prove an analogous generalization to $F$-rational pairs.
\end{rem}

We now consider another variant of $F$-rational pairs corresponding to the pure rationality defined in \ref{purelyRationalDefinition}.
\begin{defn}
\label{DivisoriallyFRationalDefinition}
Let $x$ be a non-zerodivisor of a reduced ring $R$ of characteristic $p>0$ and let $\a \subseteq R$ be an ideal which is not contained in any minimal prime of $xR$. Let $t \ge 0$ be a real number.
When $R$ is local, then the triple $(R,x;\a^t)$ is said to be \textit{divisorially $F$-rational} if $I^{{\mathrm{div}*(x;\a^t)}}=I$ for every ideal $I$ generated by  a system of parameters for $R$.
When $R$ is not local, we say that $(R,x;\a^t)$ is divisorially $F$-rational if the localization $(R_{\m},\a_{\m}^t)$ is divisorially $F$-rational for every maximal ideal $\m$ of $R$.
\end{defn}

We can prove analogues of Proposition \ref{basic}, Lemma \ref{local characterization} and Proposition \ref{basic2} for divisorial $F$-rationality.

\begin{prop}
\label{ComparisonOfFDivisorialRationalsAndFRegulars}
Let $x$ is a non-zerodivisor of a reduced ring $R$ of characteristic $p>0$ and let $\a \subseteq R$ be an ideal which is not contained in any minimal prime of $xR$. Let $t \ge 0$ be a real number.
\begin{enumerate}
\item If $(R,x;\a^t)$ is divisorially $F$-rational, then $(R,x^{1-\epsilon}\a^t)$ is $F$-rational for all $1 \ge \epsilon >0$; in particular, $(R,\a^t)$ is $F$-rational.
\item
Assume in addition that $R$ is locally excellent.
If $(R,x;\a^t)$ is divisorially $F$-regular, then it is divisorially $F$-rational.
If $R$ is $F$-finite Gorenstein, then the converse also holds true.
\end{enumerate}
\end{prop}
\begin{proof}
(2) follows from the combination of Lemma \ref{divisorial characterization} and an argument similar to the proof of Proposition \ref{basic2} (1). So, we will prove only (1).
Without loss of generality we may assume that $R$ is local.
Let $I \subseteq R$ be an ideal generated by a system of parameters for $R$ and let $z \in I^{*x^{1-\epsilon}\a^t}$.
By definition, there exists $c \in R^{\circ}$ such that $cx^{\lceil (1-\epsilon)q \rceil}\a^{\lceil tq \rceil}z^q\subseteq I^{[q]}$ for all large $q=p^e$.
Then one can choose an element $d \in R^{\circ}$ which is not in any minimal prime of $xR$ such that $dx^n$ lies in the ideal $cR$ for some $n \in \N$.
Taking sufficiently large $q=p^e$ so that $n+\lceil (1-\epsilon)q \rceil \le q-1$, one has $dx^{q-1}\a^{\lceil tq \rceil}z^q \subseteq I^{[q]}$.
This implies that $z \in I^{\mathrm{div}*(x;\a^t)}=I$, because $(R,x;\a^t)$ is divisorially $F$-rational.
Thus, $(R,x^{1-\epsilon}\a^t)$ is $F$-rational.
\end{proof}

\begin{lem}\label{divisorial characterization}
Let $(R,\m)$ be a $d$-dimensional excellent reduced local ring of characteristic $p>0$ and let $t \ge 0$ be a real number.
Fix $x \in R^{\circ}$ and let $\a \subseteq R$ be an ideal which is not contained in any minimal prime of $xR$.
Then the following three conditions are equivalent to each other.
\begin{itemize}
\item[(1)] $(R,x;\a^t)$ is divisorially $F$-rational.
\item[(2)] $R$ is equidimensional and $I^{\mathrm{div}*(x;\a^t)}=I$ for some ideal $I$ generated by a system of parameters for $R$.
\item[(3)] $R$ is Cohen-Macaulay and $0^{\mathrm{div}*(x;\a^t)}_{H^d_{\m}(R)}=0$ in $H^d_{\m}(R)$.
\end{itemize}
Furthermore, if $\ba = (f)$ is principal, then the previous three conditions are equivalent to the following:
\begin{itemize}
\item[(4)]  $R$ is Cohen-Macaulay, and for each $c \in R^{\circ}$ not contained in any minimal prime of $x$, there exists $q=p^e$ such that $cf^{\lceil tq \rceil} F^e:H^d_{\m}(R) \to H^d_{\m}(R)$ is injective.
\end{itemize}
\end{lem}
\begin{proof}
The proof is essentially the same as that of Lemma \ref{local characterization}.
\end{proof}

\begin{rem}
Let the notation be as in Lemma \ref{divisorial characterization} and assume in addition that $R$ is a homomorphic image of a Gorenstein local ring.
Then one could define the \textit{divisorial test submodule} $\tau^{\mathrm{div}}(\omega_R,x;\a^t)$ associated to $(R,x;\a^t)$ to be
$$\tau^{\mathrm{div}}(\omega_R,x;\a^t)=\mathrm{Ann}_{\omega_R}(0^{\mathrm{div}*(x;\a^t)}_{H^d_{\m}(R)}) \subseteq \omega_R.$$
This is a characteristic $p$ analogue of the adjoint submodule (see Definition \ref{adjointSubmoduleDefinition}).
By the above lemma, $(R,x;\a^t)$ is divisorially $F$-rational if and only if $\tau^{\mathrm{div}}(\omega_R,x;\a^t)=\omega_R$ and $R$ is Cohen-Macaulay.
\end{rem}

\section{Geometric Properties}
\label{SectionGeometricProperties}

In fixed prime characteristic, $F$-rational pairs satisfy several nice properties analogous to those of rational pairs.
\begin{prop}
\label{SummandTheoremForPositiveCharacteristic}
Let $R \hookrightarrow S$ be a pure finite local homomorphism of local domains of characteristic $p>0$.
Let $\a$ be a nonzero ideal of $R$ and $t \ge 0$ be a real number.
If $(S, (\a S)^t)$ is $F$-rational, then so is $(R, \a^t)$.
\end{prop}
\begin{proof}
Let $I \subseteq R$ be an ideal generated by a system of parameters for $R$.
Then it is easy to check that $I^{*\a^t}S \subseteq (I S)^{*(\a S)^t} $.
Since $I S$ is generated by a system of parameters for $S$, by assumption, $(I S)^{*(\a S)^t}=IS$.
Thus,
$$I^{*\a^t}=I^{*\a^t}S \cap R \subseteq (I S)^{*(\a S)^t} \cap R=I,$$ that is,  $(R, \a^t)$ is $F$-rational.
\end{proof}

\begin{rem}
Suppose $R$ and $S$ are domains, $\a$ is an ideal of $R$ and $R$ is a direct summand of $S$. If $(S,(\a S)^t)$ is strongly $F$-regular, then $(R,\a^t)$ is also strongly $F$-regular, in particular, $F$-rational.
However, even if $(S,(\a S)^t)$ is $F$-rational, $(R,\a^t)$ is not necessarily $F$-rational in general (see \cite{WatanabeFRationalityOfCertainReesAndCounterExamplesToBoutot} and \cite{HaraWatanabeYoshidaFRationalityOfRees} for counterexamples).
The reader should compare this with Corollary \ref{rationalSummands}.
\end{rem}

\begin{prop}\label{inversion}
Let $(R,\m)$ be an excellent reduced local ring of characteristic $p>0$ and let $x \in \m$ be a non-zerodivisor of $R$, and denote $S:=R/xR$.
Let $\a \subseteq R$ be an ideal which is not contained in any minimal prime of $xR$ and let $t \ge 0$ be a real number.
Then $(S,(\a S)^t)$ is $F$-rational if and only if $(R,x;\a^t)$ is divisorially $F$-rational.
\end{prop}
\begin{proof}
First assume that $(S,(\a S)^t)$ is $F$-rational. Note that both $S$ and $R$ are normal and Cohen-Macaulay by Proposition \ref{basic}.
We choose elements $y_1, \dots, y_{d-1}$ in $R$ such that $x, y_1, \dots, y_{d-1}$ forms a system of parameters for $R$.
Let $z \in (x, y_1, \dots, y_{d-1})^{\mathrm{div}*(x;\a^t)}$.
Then there exists $c \in R \setminus xR$ such that $c\a^{\lceil tq \rceil}x^{q-1}z^q \subseteq (x^q, y_1^q, \dots, y_{d-1}^q)$ for all large $q=p^e$.
Since $x, y_1, \dots, y_{d-1}$ is an $R$-regular sequence,
one has $c\a^{\lceil tq \rceil}z^q \in (x, y_1^q, \dots, y_{d-1}^q)$.
This implies that $\overline{c}(\a S)^{\lceil tq \rceil} \overline{z}^q \in (\overline{y_1}^q, \dots, \overline{y_{d-1}}^q)$ where $\overline{c},  \overline{z}, \overline{y_1}, \dots, \overline{y_{d-1}}$ are the images of $c, z, y_1, \dots, y_{d-1}$ in $S$, respectively.
Since $\overline{c} \in S \setminus \{0\}=S^{\circ}$,
$$\overline{z} \in (\overline{y_1}, \dots, \overline{y_{d-1}})^{*(\a S)^t}=(\overline{y_1}, \dots, \overline{y_{d-1}}).$$
Thus, $z$ lies in $(x, y_1, \dots, y_{d-1})$ and $(R,x;\a^t)$ is divisorially $F$-rational by Lemma \ref{divisorial characterization}.
The converse argument just reverses this.
\end{proof}

As a corollary of Proposition \ref{inversion}, we obtain the correspondence between pure rationality and divisorial $F$-rationality.
Let $R$ be an algebra essentially of finite type over a field of characteristic zero and let $t \ge 0$ be a real number.
Let $x$ be a non-zerodivisor of $R$ and let $\a \subseteq R$ be an ideal which is not contained in any minimal prime of $xR$.
Then $(R,x;\a^t)$ is said to be of \textit{divisorially $F$-rational} if reduction to characteristic $p \gg 0$ of $(R,x;\a^t)$ is divisorially $F$-rational (see the paragraph before Definition \ref{reduction} for the meaning of ``reduction to characteristic $p \gg 0$").
\begin{cor}
Let the notation be as above and assume in addition that $R$ is a normal local ring essentially of finite type over a field of characteristic zero.
Then $(\Spec R, \mathrm{div}(x), \a^t)$ has purely rational singularities if and only if $(R,x, a^t)$ is of divisorially $F$-rational type.
\end{cor}
\begin{proof}
It follows from the combination of Theorem \ref{correspondence}, Proposition \ref{inversion} and Theorem \ref{purelyRationalInversionOfAdjunction}.
\end{proof}

Let $R$ be a reduced ring of characteristic $p > 0$.  Suppose that $f \in R^{\circ}$ and 
that $t \ge 0$ be a real number. Let $N \subseteq M$ be $R$-modules.
Then the \textit{$\a^t$-Frobenius closure} $N^{F f^t}_M$ of $N$ in $M$ is defined to be the submodule of $M$ consisting of all elements $z \in M$ for which
$$f^{\lfloor t(q-1) \rfloor}z^q \subseteq N^{[q]}_M$$
for all large $q = p^e$. The $f^t$-Frobenius closure of an ideal $I$ of $R$ is simply defined by $I^{F f^t}:=I_R^{F f^t}$.

\begin{defn}
\label{DefinitionFInjectiveThreshold}
Let $R$ be a reduced Cohen-Macaulay ring of characteristic $p>0$ with an element $f \in R^{\circ}$ and let $t \ge 0$ be a real number.
\renewcommand{\labelenumi}{(\roman{enumi})}
\begin{enumerate}
\item
When $(R,\m)$ is a $d$-dimensional local ring and $F^e:H^d_{\m}(R) \to H^d_{\m}(R)$ denotes the induced $e$-times iterated Frobenius map on $H^d_{\m}(R)$, the pair $(R,f^t)$ is said to be \textit{$F$-injective} if for all large $q=p^e$, $f^{\lfloor t(q-1) \rfloor} F^e:H^d_{\m}(R) \to H^d_{\m}(R)$ is injective.
When $R$ is not local, we say that $(R,f^t)$ is $F$-injective if the localization $(R_{\m}, f_{\m}^t)$ is $F$-injective for every maximal ideal $\m$ of $R$.
\item Suppose that $R$ is $F$-rational. Then we define the $F$-injective threshold $\mathrm{fit}(f)$ of $f$ to be
$$\mathrm{fit}(f)=\sup\{t \in \R_{\ge 0} \mid (R,f^t) \textup{ is $F$-rational}\}.$$
\end{enumerate}
\end{defn}

We briefly study the properties of $F$-injective pairs which are needed in subsequent propositions.

\begin{lem}\label{Finjective characterization}
Let the notation be as above and assume in addition that $(R,\m)$ is a $d$-dimensional reduced Cohen-Macaulay local ring.
Then the following three conditions are equivalent to each other.
\begin{enumerate}
\item $(R,f^t)$ is $F$-injective.
\item $I^{F f^t}=I$ for some ideal $I$ generated by a system of parameters for $R$.
\item $0^{F f^t}_{H^d_{\m}(R)}=0$ in $H^d_{\m}(R)$.
\end{enumerate}
\end{lem}
\begin{proof}
The equivalence of (2) and (3) is essentially the same as in Lemma \ref{local characterization}.  Furthermore, it is clear that (1) implies (3).  To see that (3) implies (1), suppose that $0^{F f^t}_{H^d_{\m}(R)}=0$.

For each non-zero $z \in H^d_m(R)$, there exists infinitely many $e \geq 0$ such that
\[
f^{\lfloor t(p^e - 1) \rfloor} z^{p^e} \neq 0.
\]
We we show that the same statement is true for \emph{all} $e \geq 0$.

Suppose that for some $e'$, we have $f^{\lfloor t(p^{e'} - 1) \rfloor} z^{p^{e'}} = 0$.  For each $b > 0$, this implies that
\[
f^{p^b \lfloor t(p^{e'} - 1) \rfloor} z^{p^{e'+b}} = 0.
\]
But then $p^b \lfloor t(p^{e'} - 1) \rfloor \leq t(p^{e'}p^d - p^b) \leq t(p^{e'}p^b - 1)$.  The far left-side is an integer, thus we get that
$p^b \lfloor t(p^{e'} - 1)\rfloor \leq \lfloor t(p^{e'}p^b - 1)\rfloor$.
But then we also see that
\[
f^{\lfloor t(p^{e'}p^b - 1) \rfloor} z^{p^{e'+b}} = 0,
\]
a contradiction.  Therefore, we see that for each non-zero $z \in H^d_m(R)$ and for every $e > 0$, we have that
\[
f^{\lfloor t(p^e - 1) \rfloor} z^{p^e} \neq 0.
\]
\end{proof}

If one wishes to work with $F$-injective pairs $(R, \ba^t)$ when $\ba$ is not a principal ideal, the authors believe that the criterion given in  Lemma \ref{Finjective characterization}(3) is a better behaved notion than the definition we gave above in Definition \ref{DefinitionFInjectiveThreshold}(1).  We will restrict ourselves to pairs $(R, f^t)$ in this paper however.


\begin{lem}\label{F-inj lem}
Let the notation be as in Lemma \ref{Finjective characterization}
\begin{enumerate}
\item If $(R,f^t)$ is $F$-rational, then it is $F$-injective.
\item Suppose that $R$ is an excellent $F$-rational local ring. Then $(R,f^t)$ is $F$-injective if and only if $t \le \mathrm{fit}(f)$.
\end{enumerate}
\end{lem}
\begin{proof}
(1) is immediate from definition.
So, we will prove only (2). 

First we will show that
$$\mathrm{fit}(f)=\sup\{t \in \R_{\ge 0} \mid (R,f^t) \textup{ is $F$-injective}\}.$$
To check this, it is enough to show that if $(R,f^t)$ is $F$-injective for some $t>0$, then $(R,f^{t-\epsilon})$ is $F$-rational for all $t \ge \epsilon >0$.
By Lemma \ref{test exist}, the unit $1$ is a parameter $f^{t-\epsilon}$-test element.
Since $\lceil (t-\epsilon)q \rceil \le \lfloor t(q-1) \rfloor$ for all large $q=p^e$, we have the inclusion $I^{*f^{t-\epsilon}} \subseteq I^{F f^t}$ for every ideal $I$ generated by a system of parameters for $R$.
Thus, the F-injectivity of $(R,f^t)$ implies the F-rationality of $(R,f^{t-\epsilon})$.

To complete the proof of (2), 
it only remains to show that $(R,f^{\mathrm{fit}(f)})$ is $F$-injective.
Let $I \subseteq R$ be an ideal generated by a system of parameters for $R$.
Let
$$\nu(p^e):=\max\{r\in\N\vert f^rz^{p^e}\not\subseteq I^{[p^e
]} \textup{ for all $z \in R \setminus I$}\}.$$
Since $R$ is $F$-injective, the invariant $\nu(p^e)$ is well-defined.
\begin{cl}[cf.~\cite{MustataTakagiWatanabeFThresholdsAndBernsteinSato}]
$$\mathrm{fit}(f)=\lim_{e \to \infty} \frac{\nu(p^e)}{p^e}=\inf_{e} \frac{\nu(p^e)+1}{p^e}.$$
\end{cl}
\begin{proof}[Proof of Claim]
Since $f$ is a principal ideal, if $f^{\nu(q)+1}z^q$ lies in $I^{[q]}$, then $f^{p(\nu(q)+1)}z^{pq}$ lies in $I^{[pq]}$.
Thus, $({\nu(q)+1})/{q} \ge ({\nu(pq)+1})/{pq}$, that is,
$$\lim_{e \to \infty} \frac{\nu(p^e)}{p^e}=\inf_{e} \frac{\nu(p^e)+1}{p^e}.$$
Since $(R,f^{\mathrm{fit}(f)-\epsilon})$ is $F$-injective and $(R,f^{\mathrm{fit}(f)+\epsilon})$ is never $F$-injective for all $1 \ge \epsilon >0$, by (2), one has
$\lfloor (\mathrm{fit}(f)-\epsilon)(q-1) \rfloor \le \nu(q) < \lfloor (\mathrm{fit}(f)+\epsilon)(q-1) \rfloor$
for infinitely many $q=p^e$. This implies that
$$\mathrm{fit}(f)-\epsilon \le \lim_{e \to \infty} \frac{\nu(p^e)}{p^e} \le \mathrm{fit}(f)+\epsilon.$$
Since $\epsilon$ can take arbitrary small values, we obtain the assertion.
\end{proof}

By the above claim, $\lfloor \mathrm{fit}(f)(q-1) \rfloor \le  \nu(q)$ for every $q=p^e$, which means that $I^{Ff^{\mathrm{fit}(f)}}=I$.
\end{proof}

\begin{thm}
\label{FInjectiveInversion}
Let $(R,\m)$ be an excellent reduced local ring of characteristic $p>0$ and $x \in \m$ be a non-zerodivisor of $R$.
If $(R,x^{1-\epsilon})$ is $F$-rational for all sufficiently small $1 \gg \epsilon>0$, then $R/xR$ is Cohen-Macaulay and $F$-injective $($that is, the pair $(R/xR, (R/xR)^1)$ is $F$-injective$)$.
When the localized ring $R_x$ is $F$-rational, the converse implication also holds true.
\end{thm}
\begin{proof}
Without loss of generality, we may assume that $R$ is Cohen-Macaulay.

\begin{cl}
$(R,x)$ is $F$-injective if and only if $R/xR$ is $F$-injective.
\end{cl}
\begin{proof}[Proof of Claim]
We choose elements $y_1, \dots, y_{d-1}$ in $R$ such that $x, y_1, \dots, y_{d-1}$ is a system of parameters for $R$. An element $z \in R$ lies in $(x, y_1, \dots, y_{d-1})^{F x}$ if and only if $z^q \in (x, y_1^q, \dots, y_{d-1}^q)$ for all large $q=p^e$, because $x, y_1, \dots, y_{d-1}$ is an $R$-regular sequence.
This is equivalent to saying that $\overline{z}^q \in (\overline{y_1}^q, \dots, \overline{y_{d-1}}^q)$ for all large $q=p^e$, that is, $\overline{z} \in (\overline{y_1}, \dots, \overline{y_{d-1}})^F$, where $\overline{z}, \overline{y_1}, \dots, \overline{y_{d-1}}$ are the images of $z, y_1, \dots, y_{d-1}$ in $R/xR$
, respectively.
Thus, 
we obtain the assertion.
\end{proof}

If $(R,x^{1-\epsilon})$ is $F$-rational for all sufficiently small $1 \gg \epsilon>0$, then by Lemma \ref{F-inj lem} (2), 
$(R,x)$ is $F$-injective.
To complete the proof of this theorem, by the above claim, it only remains to show that if $(R,x)$ is $F$-injective and $R_x$ is $F$-rational, then  $(R,x^{1-\epsilon})$ is $F$-rational for all $1 \ge \epsilon >0$.
Since $R_x$ is $F$-rational, by Lemma \ref{test exist}, some power $x^n$ of $x$ is a parameter $x$-test element.
Since $\lceil (1-\epsilon)q \rceil+n \le q-1$ for all large $q=p^e$, if $I \subseteq R$ is an ideal generated by a system of parameters for $R$ and $z \in I^{F x}$, then $x^n x^{\lceil (1-\epsilon)q \rceil} z^q \in I^{[q]}$ for all large $q=p^e$, that is, $z \in I^{*x^{1-\epsilon}}$.  Thus, the F-injectivity of $(R,x)$ implies the F-rationality of $(R,x^{1-\epsilon})$.
\end{proof}

\begin{eg}
\label{VeroneseExamples}
Consider the $r$-th Veronese subring $R=S^{(r)}$ of the d-dimensional formal power series ring $S=k[[x_1, \dots, x_d]]$ over a perfect field $k$ of characteristic $p>0$. It is well-known that $R$ is strongly $F$-regular.
By \cite[Example 2.4 (ii)]{TakagiWatanabeFPureThresh}, the $F$-pure threshold $\mathrm{fpt}(\m)$ of the maximal ideal $\m$ of $R$ is equal to $d/r$, that is, $(R,\m^t)$ is strongly $F$-regular if and only if $t<d/r$. We will show that the $F$-injective threshold $\mathrm{fit}(\m)$ of $\m$ is equal to $\lceil d/r \rceil$.

Let $I=(x_1^r, x_2^r, \dots, x_d^r)$.  Then $(R,\m^t)$ is $F$-rational if and only if $I^{*\m^t}$ contains none of the monomials $x_1^{i_1} \dots x_d^{i_d}$ in $R$ with $r-2 \le i_j \le r-1$ for all $j=1, \dots, d$.
Put $n=\lceil d/r \rceil$ and $z=x_1^{r-2} \cdots x_{rn-d}^{r-2}x_{rn-d+1}^{r-1} \cdots x_d^{r-1} \in R$.
Since $x_1^{2q-2} \cdots x_{rn-d}^{2q-2} x_{rn-d+1}^{q-1} \cdots x_d^{q-1}$ is in $\m^{n(q-1)} \subseteq \m^{\lceil (n-\epsilon)q \rceil}$ for all large $q=p^e$, $z^q \m^{\lceil (n-\epsilon)q \rceil}$ is not contained in $I^{[q]}$. Thus, $z$ does not belong to $I^{*\m^{n-\epsilon}}$ (here, $1$ is an $\m^{n-\epsilon}$-test element by Lemma \ref{completion}). Similarly, we can show that $I^{*\m^{n-\epsilon}}$ contains none of the monomials $x_1^{i_1} \dots x_d^{i_d}$ in $R$ with $r-2 \le i_j \le r-1$ for all $j=1, \dots, d$.
This means $\mathrm{fit}(\m) \ge n=\lceil d/r \rceil$. We leave it for the reader to check that $\mathrm{fit}(\m) \le \lceil d/r \rceil$.
\end{eg}

We conclude this section with a proof of a special case of the discreteness and rationality of $F$-injective thresholds.
More generally, we introduce a new invariant which is a generalization of $F$-injective thresholds and study its properties.

\begin{defn}
\label{DefinitionJumpingExponentForTestModules}
Let $R$ be a reduced local ring of characteristic $p>0$ which is a homomorphic image of a Gorenstein local ring and let $\a \subseteq R$ be an ideal such that $\a \cap R^{\circ} \ne \emptyset$.
We say that a real number $t>0$ is a \textit{jumping exponent for generalized parameter test submodules $\tau(\omega_R, \a^*)$} if $\tau(\omega_R, \a^{t}) \subsetneq \tau(\omega_R, \a^{t-\epsilon})$ for all $\epsilon >0$.
\end{defn}

If $R$ is excellent and $F$-rational, then by Remark \ref{testModuleDefinition}, the smallest jumping exponent for the generalized parameter test submodules $\tau(\omega_R, \a^*)$ is the $F$-injective threshold $\mathrm{fit}(\a)$ of $\a$.

\begin{lem}\label{conti}
Let $(R,\m)$ be a complete local domain of characteristic $p>0$ and $\a \subseteq R$ be an ideal such that $\a \cap R^{\circ} \ne \emptyset$.
\renewcommand{\labelenumi}{$(\arabic{enumi})$}
\begin{enumerate}
\item
For every nonnegative real number $t$, there exists $\epsilon >0$ such that
$$\tau(\omega_R, \a^t)=\tau(\omega_R, \a^{t'})$$ for all $t' \in [t, t+\epsilon)$.
\item
If $\alpha$ is a jumping exponent for generalized parameter test submodules $\tau(\omega_R, \a^*)$, then so is $p \alpha$.
\item
If $\a$ is generated by $m$ elements, then for every $t \ge m$, one has
$$\tau(\omega_R, \a^t)=\tau(\omega_R, \a^{t-1}) \a.$$
\end{enumerate}
\end{lem}

\begin{proof}
(1)
Let $c \in R^{\circ}$ be a parameter $\a^s$-test element for every $s \ge 0$, and fix any $d \in \a \cap R^{\circ}$.
Then $cd$ is also a parameter $\a^s$-test element for every $s \ge 0$.
Denote by $N_e$ the submodule of $H^{d}_{\m}(R)$ consisting of all elements $\xi \in H^{d}_{\m}(R)$ such that $cd \a^{\lceil tp^e \rceil}\xi^{p^e}=0$ in $H^{d}_{\m}(R)$.
By definition, one can see that $0^{*\a^t}_{H^d_{\m}(R)}=\bigcap_{e \in \N}N_e$.
Since $H^d_{\m}(R)$ is an Artinian $R$-module, there exists an integer $m$ such that $0^{*\a^t}_{H^d_{\m}(R)}=\bigcap_{e=0}^mN_e$.
Put $\epsilon=1/p^m$ and we will prove that $\tau(\a^t)=\tau(\a^{t+\epsilon})$.
Let $\xi \in 0_{H^{d}_{\m}(R)}^{*\a^{t+\epsilon}}$.
Since $c$ is a parameter $\a^{t+\epsilon}$-test element, $c\a^{\lceil (t+\epsilon)q \rceil}\xi^q=0$ for all $q=p^e$.
In particular, $c\a^{\lceil tp^e \rceil+1}\xi^{p^e}=0$ for all $e=0, \dots, m$.
Since $d$ is in $\a$,  $\xi$ lies in $\bigcap_{e=0}^m N_e=0^{*\a^t}_{H^d_{\m}(R)}$.

(2)
Let $c \in R^{\circ}$ be a parameter $\a^t$-test element for every $t \ge 0$, and fix any $\epsilon>0$.
Since $\alpha$ is a jumping exponent for generalized parameter test submodules $\tau(\omega_R, \a^*)$, there exists $\xi \in 0_{H^{d}_{\m}(R)}^{*\a^{\alpha}}$ which is not contained in $0_{H^{d}_{\m}(R)}^{*\a^{\alpha-\epsilon}}$.
This means that $c\a^{\lceil \alpha q \rceil}\xi^q=0$ in $H^{d}_{\m}(R)$ for all $q=p^e$, but $c\a^{\lceil (\alpha-\epsilon)q \rceil}\xi^q \ne 0$ in $H^{d}_{\m}(R)$ for infinitely many $q=p^e$.
Put $\eta=\xi^p \in H^{d}_{\m}(R)$.
Then $c\a^{\lceil p\alpha q \rceil}\eta^q=0$ in $H^{d}_{\m}(R)$ for all $q=p^e$, but $c\a^{\lceil p(\alpha-\epsilon)q \rceil}\eta^q \ne 0$ for infinitely many $q=p^e$.
This implies that $\eta$ belongs to $0_{H^{d}_{\m}(R)}^{*\a^{p\alpha}}$ but not to $0_{H^{d}_{\m}(R)}^{*\a^{p(\alpha-\epsilon)}}$.
Thus, by Matlis duality, $\tau(\omega_R, \a^{p\alpha}) \subsetneq \tau(\omega_R, \a^{p(\alpha-\epsilon)})$.

(3)
By the proof of \cite[Theorem 4.1]{HaraTakagiOnAGeneralizationOfTestIdeals}, $0^{*\a^t}_{H^d_{\m}(R)} = (0^{*\a^{t-1}}_{H^d_{\m}(R)} \colon \a)_{H^d_{\m}(R)}$ for every real number $t \ge m$.
Since $\Ann_{H^d_{\m}(R)} (\tau(\omega_R, \a^{t-1})\a)$ is equal to $(0^{*\a^{t-1}}_{H^d_{\m}(R)}  \colon \a)_{H^d_{\m}(R)}  = 0^{*\a^t}_{H^d_{\m}(R)} $,  by Matlis duality, one has
$$\tau(\omega_R, \a^t)
            =  \Ann_{\omega_R}(0^{*\a^t}_{H^d_{\m}(R)})
            = \Ann_{\omega_R}(\Ann_{H^d_{\m}(R)} (\tau(\omega_R, \a^{t-1})\a)) = \tau(\omega_R, \a^{t-1})\a.$$
\end{proof}

\begin{thm}\label{exponent}
Let $(R,\m)$ be an excellent $F$-rational local ring of characteristic $p>0$ which is a homomorphic image of a Gorenstein local ring, and fix $g \in R^{\circ}$.
Let $\alpha, \beta>0$ be integers and write $\gamma=\alpha/(p^{\beta}-1)$.
Then there exists $c \in (0, \gamma)$ for which $\tau(\omega_R, g^t)=\tau(\omega_R, g^c)$ for all $t \in [c, \gamma)$.
\end{thm}
\begin{proof}
For any integers $m, n>0$, we denote by $N_{m,n}$ the submodule of $H^{d}_{\m}(R)$ consisting of all elements $\xi \in H^{d}_{\m}(R)$ such that $g^{m}\xi^{p^n}=0$ in $H^{d}_{\m}(R)$.

\begin{cl}
$\tau(\omega_R, g^{m/p^n})=\Ann_{\omega_R} (N_{m,n})$.
\end{cl}
\begin{proof}[Proof of Claim]
First note that by Lemma \ref{test exist}, the unit $1$ is a parameter  $g^t$-test element for every $t \ge 0$.
Then one can see that $0^{*g^{m/p^n}}_{H^d_{\m}(R)}=\bigcap_{e \ge n}N_{mp^{e-n},e}$.
Now it suffices to show that $N_{k,e} \subseteq N_{kp,e+1}$ for all integers $k, e>0$, but it is obvious.
\end{proof}

Fix the $R[\theta; f^\beta]$-module structure on $H^d_{\m}(R)$ given by $\theta \xi=g^{\alpha}\xi^{p^{\beta}}$ for all $\xi \in H^d_{\m}(R)$.
Then $N_{\alpha(1+p^{\beta}+\dots+p^{\beta(s-1)}), s\beta}$ coincides with the kernel of $\theta^s$ as an  $R[\theta; f^{\beta}]$-module.
Thus, $\{N_{\alpha(1+p^{\beta}+\dots+p^{\beta(s-1)}), s\beta}\}_{s \ge 1}$ forms an ascending chain of $R[\theta; f^{\beta}]$-modules, and by Hartshorne-Speiser-Lyubeznik's theorem (\cite[Proposition 4.4]{LyubeznikFModulesApplicationsToLocalCohomology}) it stabilizes at some $s=\nu$.
For all $s \ge 1$, by the above claim,
$$\tau\left(\omega_R, g^{\frac{\alpha(1+p^{\beta}+\dots+p^{\beta(s-1)})}{p^{s\beta}}}\right)=\Ann_{\omega_R}(N_{\alpha(1+p^{\beta}+\dots+p^{\beta(s-1)}), s\beta}).$$
Thus, the stabilization of $\{N_{\alpha(1+p^{\beta}+\dots+p^{\beta(s-1)}), s\beta}\}_{s \ge 1}$ implies that of a family of generalized parameter test submodules $\{\tau(\omega_R, g^{\alpha(1+p^{\beta}+\dots+p^{\beta(s-1)})/{p^{s\beta}}})\}_{s \ge 1}$ for $s \ge \nu$.
Since $$\frac{\alpha(1+p^{\beta}+\dots+p^{\beta(s-1)})}{p^{s\beta}}=\frac{\alpha}{p^{s\beta}}\frac{p^{s\beta}-1}{p^{\beta}-1}$$
is an increasing sequence which converges to $\gamma$ as $s$ goes to the infinity,
we may take $c=\alpha(1+p^{\beta}+\dots+p^{\beta(\nu-1)})/{p^{\nu\beta}}$.
\end{proof}

\begin{cor}\label{accumulation}
Let $(R,\m)$ be a complete $F$-rational local ring of characteristic $p>0$, and fix $g \in R^{\circ}$. Then the set of jumping exponents for generalized parameter test submodules $\tau(\omega_R, g^*)$ cannot have a rational accumulation point.
\end{cor}
\begin{proof}
Assume to the contrary that the set of jumping exponents for generalized parameter test submodules $\tau(\omega_R, g^*)$ have a rational accumulation point $\gamma$.
Then there exists a sequence $\{c_n\}_{n \ge 1}$ of jumping exponents converging to $\gamma$.
When we write $\gamma$ in the form of $\frac{\alpha}{p^d(p^{\beta}-1)}$, by Lemma \ref{conti} (2), $\{p^d c_n\}_{n \ge 1}$ is a sequence of jumping exponents again and is converging to $\frac{\alpha}{p^{\beta}-1}$. This contradicts Theorem \ref{exponent}.
\end{proof}

\begin{cor}\label{rational}
Let $(R,\m)$ be a complete $F$-rational local ring of characteristic $p>0$, and fix $g \in R^{\circ}$. Then jumping exponents for generalized parameter test submodules $\tau(\omega_R, g^*)$ are rational and have no  accumulation points.
\end{cor}
\begin{proof}
Applying Lemma \ref{conti} and Corollary \ref{accumulation} to \cite[Proposition 4.2]{KatzmanLyubeznikZhangOnDiscretenessAndRationality}, we obtain the assertion.
\end{proof}

\providecommand{\bysame}{\leavevmode\hbox to3em{\hrulefill}\thinspace}
\providecommand{\MR}{\relax\ifhmode\unskip\space\fi MR}
\providecommand{\MRhref}[2]{
  \href{http://www.ams.org/mathscinet-getitem?mr=#1}{#2}
}
\providecommand{\href}[2]{#2}

\end{document}